\documentclass[EJP,preprint]{ejpecp}
\usepackage{booktabs,array,caption,xcolor,xurl}
\usepackage[shortlabels]{enumitem}
\hypersetup{colorlinks=true,linkcolor=blue!45!black,citecolor=blue!45!black,urlcolor=blue!45!black}

\newcommand{\E}{\mathbb{E}}
\newcommand{\Prob}{\mathbb{P}}
\newcommand{\Var}{\operatorname{Var}}

\SHORTTITLE{The bias interval of the elephant random walk}
\TITLE{The bias interval of the elephant random walk}
\hypersetup{pdftitle={The bias interval of the elephant random walk},
pdfauthor={Arjun Pemmasani}}
\AUTHORS{Arjun Pemmasani\footnote{\raggedright Department of Mathematics, Harvey Mudd College,
Claremont, CA 91711, USA. \mbox{\EMAIL{apemmasani@hmc.edu}}\newline
ORCID \href{https://orcid.org/0009-0001-2776-6298}{0009-0001-2776-6298}.}}
\AMSSUBJ{60G50}
\AMSSUBJSECONDARY{60F05; 60C05}
\ABSTRACT{Gu\'erin, Laulin, Raschel, and Simon conjectured that the
initial biases giving the superdiffusive elephant random walk a
bimodal limit law form an interval centered at one half, and they
described the problem as challenging. We prove this. For every
superdiffusive memory the limit density has two local maxima exactly
when the bias lies in a proper open interval, and outside it,
including its endpoints, the density has a unique local maximum. The
interval is nonempty exactly above a single memory threshold,
numerically about $0.6688$. Near this threshold $a_*$ the distance
between the two modes is of order $(a-a_*)^{1/2}$, and the width of
the interval is of order $(a-a_*)^{3/2}$. The bias boundary is
continuous and tends to zero at least quadratically fast as memory
tends to one. The proof bounds the sign changes of the first and
second differences of the finite walk, including those created at the
boundary, and passes them to the limit.}
\begin{document}
\raggedbottom
\section{Introduction}\label{sec:intro}\label{sec:notation}

\noindent\begin{minipage}{\linewidth}
\begin{quote}\small\itshape
An elephant and an ant are playing hide and seek. The elephant counts while
the ant hides among many identical anthills, then opens its eyes and knows
at once which one the ant is in. The ant's shoes are outside the door.
\end{quote}
\hfill\textup{(Anonymous)}
\medskip
\end{minipage}

The elephant random walk repeats or reverses one of its own earlier
steps, chosen uniformly. If its first step is $+1$ with probability
$q$, for which $q$ does its limit law have two peaks? When the memory is
strong enough the walk is superdiffusive, and its rescaled position
converges to a random limit that is not
Gaussian~\cite{dasilva2013}. Given that the first step is $+1$, this
limit has a density with one peak~\cite{glrs2025}. By the symmetry of the walk, with initial bias
$q$ the limit density is a mixture of this density, with weight $q$,
and of its mirror image, with weight $1-q$ (Figure~\ref{fig:phase}).
When the memory parameter $a$ is close to one the one-peak density lies
far from zero, so the fair mixture has two peaks, and when $a$ is close
to $1/2$ the density is spread out and the two peaks merge. But a
mixture of a one-peak density and its mirror image can have more than
two peaks, and the weights that give two peaks do not have to form an
interval. The weights giving two peaks
are known for normal
mixtures~\cite{eisenberger1964,robertsonfryer1969}, but here the
density is not Gaussian. Gu\'erin, Laulin, Raschel, and Simon
conjecture that the biases $q$ with two peaks form an open interval
centered at $1/2$, and they describe this problem as
challenging~\cite[Remark~2.1(b)]{glrs2025}. In a growing-network
model, Mohandas et al.~\cite[Section~V and Figure~5]{mohandas2025}
found the fair-start transition numerically. They extrapolated
finite-$N$ central curvatures and found it at about $0.836$ in their
parameter, the probability $(1+a)/2$ that a step copies an earlier
one. We prove the conjecture, we show that the interval is nonempty
exactly above one threshold $a_*\approx0.6688$, which is
$(1+a_*)/2\approx0.8344$ in that parameter, and we describe how the
interval opens.

The elephant walk~\cite{schutz2004} has steps $X_i\in\{+1,-1\}$.
Given $X_1,\dots,X_t$, it picks $i$ uniformly from $\{1,\dots,t\}$.
Then $X_{t+1}=X_i$ with probability $1-\varepsilon$, and
$X_{t+1}=-X_i$ otherwise. These choices are independent of each other
and of $X_1$. Put
\[
 A_N=\#\{i\le N:X_i=+1\},\qquad S_N=2A_N-N,\qquad
 P_N(k)=\Prob(A_N=k).
\]
The memory is $a=1-2\varepsilon$ and the initial bias is
$q=\Prob(X_1=1)$. We call $q=1/2$ a fair start.

When $a\in(1/2,1)$ the ratio $S_N/N^a$ converges almost surely
\cite{baurbert2016,bercu2018}. Let $h_a$ be the density of the limit
given $X_1=1$. With initial bias $q$ the limit density is
\[
 g_{a,q}(u)=q h_a(u)+(1-q)h_a(-u),\qquad g_a=g_{a,1/2}.
\]

\subsection{The bias interval and its threshold}\label{sec:results}

\begin{theorem}[the bias interval]\label{thm:main-limit}
For every $a\in(1/2,1)$ there is $b_a\in(0,1/2]$ such that
$g_{a,q}$ has exactly two local maxima if and only if
$b_a<q<1-b_a$. Outside this open interval, including its endpoints,
the density is unimodal with a unique local maximum.
\end{theorem}

\begin{figure}[t]
\centering
\includegraphics[width=\textwidth]{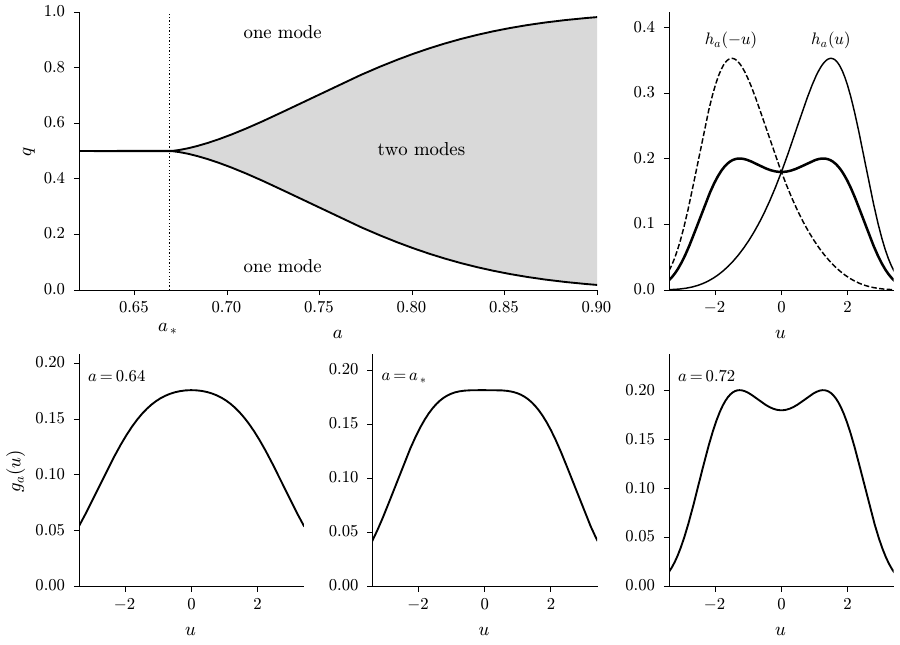}
\caption{The bias interval and the fair-start transition. Upper left:
$b_a$ and $1-b_a$ from Theorem~\ref{thm:bias-continuity}, computed
numerically, with the biases that give two modes shaded. Upper right: the density $h_a$ of the limit given
$X_1=1$ (solid), its mirror image (dashed), and their fair mixture
$g_a$ (thick), for $a=0.72$. Lower row: $g_a$ for $a=0.64$, $a=a_*$, and $a=0.72$. The top
at $a_*$ is quartic by Lemma~\ref{lem:analytic-quartic}.}
\label{fig:phase}
\end{figure}

This is Theorems~\ref{thm:two-modes}, \ref{thm:biased-limit-modes},
and~\ref{thm:bias-continuity}. The next theorem gives the set of $a$
with $b_a<1/2$ (Figure~\ref{fig:phase}, lower row). It is
Theorem~\ref{thm:threshold}.

\begin{theorem}[the analytic threshold]\label{thm:main-analytic}
There is exactly one $a_*\in(1/2,1)$ with $g_{a_*}''(0)=0$, and it
lies in $(2/3,69/100)$. The central curvature is negative below $a_*$
and positive above it, and $\partial_a g_a''(0)|_{a=a_*}>0$.
The bias interval is empty exactly when $a\le a_*$.
\end{theorem}

The bias boundary is regular (Theorem~\ref{thm:bias-continuity},
Corollary~\ref{cor:bias-boundary-rate}, and
Proposition~\ref{prop:quadratic-rate}). The function $a\mapsto b_a$ is
continuous on $(1/2,1)$, continuously differentiable on $(a_*,1)$, and
strictly decreasing immediately above $a_*$, and
\[
 b_a\to0\quad(a\uparrow1),\qquad
 \limsup_{a\uparrow1}\frac{b_a}{(1-a)^2}\le82 .
\]
We do not know if $b_a$ is decreasing on all of $(a_*,1)$
(Section~\ref{sec:open}).

At $a_*$ the fair density has a quartic top
(Lemma~\ref{lem:analytic-quartic}), and the next theorem describes how
the two modes and the bias interval open from it. It is
Propositions~\ref{prop:critical-shape} and~\ref{prop:bias-cusp}.

\begin{theorem}[the critical exponents]\label{thm:main-critical}
There are positive constants $K_*,Q_*$ such that the positive mode
$\kappa_a$ of $g_a$ and the bias boundary satisfy
\[
 \kappa_a\sim K_*\sqrt{a-a_*},\qquad
 \tfrac12-b_a\sim Q_*(a-a_*)^{3/2}\quad(a\downarrow a_*).
\]
\end{theorem}

The constants $K_*$ and $Q_*$ are explicit in three derivatives of the
densities at $a_*$ (Section~\ref{sec:critical-density}), and
Section~\ref{sec:verify} gives their numerical values.
Theorem~\ref{thm:two-modes} also counts the inflection points of
$g_a$. All proofs are written proofs and use no computer.

\subsection{From finite signs to the limiting shape}\label{sec:architecture}

Gu\'erin, Laulin, Raschel, and Simon prove that $h_a$ has one peak by
an induction on the number of steps, for the count law given $X_1=1$.
We run the same induction for every initial bias, and also for second
differences when the start is fair. In both cases the induction fails
at the two ends of the count law. We control the values there, and
this is the first new part of the proof. The second new part is a formula for $g_a''(0)$ as
an integral with a positive weight, which we use to show that
$g_a''(0)$ changes sign only once. Passing to the limit is routine, and so is the
regularity of the densities.

\emph{Counting sign changes.} Away from the two ends, one step of the
walk replaces each difference of the count law by a positive
combination of two neighboring differences, and such a step cannot add
a sign change. This is classical variation
diminution~\cite{karlin1968,karlin1964}. At the two ends a new value
enters at each step, and it can add a sign change. So we fix the final
length $N$ and carry, on both sides of the differences at time $t$, the
end values that will enter up to time $N$
(Figure~\ref{fig:augmented-rows}). Once we have bounded the sign
changes of these end values on their own, every new end value is
already in the list, and the induction runs over the longer list with
no new sign change. In this way we show that for a fair start the
second difference changes sign at most 4 times, and that for every
bias the first difference changes sign at most 3 times. The first
bound passes to $g_a''$, and the second shows in the limit that
$g_{a,q}$ does not have two strict valleys in a row, in the sense
defined next.

\emph{The interval.} We say that a function has a strict valley if it
takes at some point a value strictly below its values at a point on
each side. Let $q\ge1/2$. On $[0,\infty)$ the density $g_{a,q}$ has no
strict valley, and $g_{a,q}(-x)\le g_{a,q}(x)$ for $x\ge0$. Hence $g_{a,q}$ has two
peaks if and only if $g_{a,q}(x)>g_{a,q}(y)$ for some $x<y\le0$. For
fixed $x,y$ this difference is affine in $q$ and is $\le0$ at $q=1$. If
it is positive at some $q$, it is positive on $[1/2,q]$. Taking the
union over $x,y$ and reflecting gives an interval $(b_a,1-b_a)$. Since
$h_a'>0$ on $(-\infty,0]$ for every $a\in(1/2,1)$, a compactness
argument gives $b_a>0$ (Theorem~\ref{thm:bias-continuity}).

\emph{The threshold.} Let $E$ be an independent exponential variable.
The moments of $E^a\mathcal L_a$, where $\mathcal L_a$ is the limit of
$S_N/N^a$ given $X_1=1$, give the transform
of~\cite[Proposition~3.2]{glrs2025}. Two integrations by parts in its
Mellin transform give $g_a''(0)$ as a positive multiple of
$\mathcal R(a)-1$, where $\mathcal R(a)$ is a weighted mean against a
positive weight. Rational bounds give the sign of $g_a''(0)$ for
$a\notin[2/3,69/100]$. On $[2/3,69/100]$ the logarithmic derivative of
the weight in $a$ is increasing in the integration variable, and a
covariance identity gives $\mathcal R'(a)>0$. Hence $g_a''(0)$ changes
sign exactly once as $a$ increases.

\emph{The critical exponents.} Both exponents can be read off from the
cubic $c(t)=C_*(a-a_*)t-B_*t^3/6$, the leading part of $g_a'(t)$ near
zero, where $C_*=\partial_ag_a''(0)|_{a=a_*}$ and
$B_*=-g_{a_*}^{(4)}(0)$ are positive. Its positive zero is
$\sqrt{6C_*(a-a_*)/B_*}$, the leading term of $\kappa_a$. Its maximum
over $t>0$ is $(2\sqrt2/3)C_*^{3/2}B_*^{-1/2}(a-a_*)^{3/2}$, and
dividing by $2H_*$ with $H_*=h_{a_*}'(0)$ gives $Q_*(a-a_*)^{3/2}$.

\emph{Where the steps are proved.} The recursion for the differences is
Lemma~\ref{lem:D-rec}, and the step away from the ends is
Lemma~\ref{lem:interleave}. Lemma~\ref{lem:level-one} bounds the end
values for second differences, and the proof of
Proposition~\ref{prop:all-bias-modes} bounds them for first differences.
The two finite bounds are Theorem~\ref{thm:five-blocks} and
Proposition~\ref{prop:all-bias-modes}, and
Lemmas~\ref{lem:derivative-limits} and~\ref{lem:sp-limit} pass them
to the limit.
Proposition~\ref{prop:dominant-half} excludes a strict valley on the
right half, and Theorem~\ref{thm:biased-limit-modes} gives the interval.
Lemma~\ref{lem:balanced-start} proves $h_a'(t)>0$ for $t<0$, and
Lemma~\ref{lem:positive-conditional-mode} proves $h_a'(0)>0$. The
formula for the curvature is Lemma~\ref{lem:curvature-weight}, and
Appendix~\ref{app:analytic-threshold} proves
Theorem~\ref{thm:threshold} from it.

\subsection{Prior work}\label{sec:prior}

The counting chain predates the elephant walk. The Drezner--Farnum
model~\cite{dreznerfarnum1993}, in the formulation of
James, James, and Qi~\cite[p.~2339]{jamesjamesqi2008}, has
next-success probability $(1-\theta)p+\theta k/t$.
With $p=1/2$ and $\theta=a$ this is the probability
$\varepsilon+ak/t$ that the elephant walk steps $+1$ after $t$ steps
of which $k$ were $+1$. Their family contains
bimodal laws. Heyde~\cite{heyde2004} gives the limit regimes at fixed
parameters, and James, James, and Qi allow $\theta$ to depend on $t$.
Baur and Bertoin~\cite{baurbert2016} prove functional limit theorems
through P\'olya-type urns. The memory of the elephant walk is a random recursive
tree~\cite{kursten2016}, and the split of this tree into two subtrees
is used in Appendix~\ref{app:modes}. Zeros and recurrence of the walk
are studied in~\cite{bertoin2022,qin2025}. See
also~\cite{coletti2017,bercu2018,kubotatakei2019}.

We use the density and moments from~\cite[Theorems~1.3 and~1.4]{glr2026}
and conditional unimodality from~\cite[Theorem~1.1]{glrs2025}. Kiss
and Vet\H{o}~\cite{kissveto2022} compute the first moments for
general steps. Transition probabilities, Gaussian rates, and moment
recurrences are in the preprints
\cite{peresqin2026,muraharatokumitsu2026,vishwakarma2026}.
Trinas and Valle~\cite{tv2026} prove a superdiffusive local limit
theorem, and show the conditional law is not log-concave. Their
difference recurrence also gives Lemma~\ref{lem:derivative-limits}.
Since their paper is a preprint, we include the proof. Trinas~\cite{trinas2026}
uses their local limit theorem for the last-return law, in which the
value $g_{a,q}(0)$ appears.

\section{Finite shape bounds}\label{sec:convex}

This section proves the finite bounds behind
Theorems~\ref{thm:main-limit} and~\ref{thm:two-modes}. For a fair
start, Theorem~\ref{thm:five-blocks} shows that the second
difference of the count law changes sign at most four times. For
every initial bias, Proposition~\ref{prop:all-bias-modes} shows
that the count law has at most two modal plateaus. The limit
transfer is Lemma~\ref{lem:derivative-limits} in
Appendix~\ref{app:modes}.

We first treat second differences under a fair start. Let
$0<\varepsilon<\frac12$, $a=1-2\varepsilon$, and $r_t(k)=\varepsilon+ak/t$. For $t\ge1$ and every integer $k$ we write
$P_t(k)=\Prob(A_t=k)$, which is $0$ unless $0\le k\le t$, and
\[
 D_t(k)=P_t(k+1)-2P_t(k)+P_t(k-1).
\]
We also put $s_t=D_t(1)=P_t(0)-2P_t(1)+P_t(2)$ for $t\ge2$. Thus $D_t$ is the second difference of the count law, and $s_t$ is
its value next to an end.

\begin{lemma}[recurrence for the second difference]\label{lem:D-rec}
For every $t\ge1$ and every integer $k$,
\begin{align}
 P_{t+1}(k)&=\bigl(1-r_t(k)\bigr)P_t(k)+r_t(k-1)\,P_t(k-1),\label{eq:P-rec}\\
 D_{t+1}(k)&=\bigl(1-r_t(k+1)\bigr)D_t(k)+r_t(k-2)\,D_t(k-1).\label{eq:D-rec}
\end{align}
Also $P_t(k)=P_t(t-k)$ and $D_t(k)=D_t(t-k)$.
\end{lemma}

\begin{proof}
Given $A_t=k$, the copied step is $+1$ with probability $k/t$, so step $t+1$
is $+1$ with probability $(1-\varepsilon)\frac kt+\varepsilon(1-\frac
kt)=r_t(k)$. This gives \eqref{eq:P-rec} for $0\le k\le t+1$. At $k=0$ and
$k=t+1$, the extra term $r_t(-1)P_t(-1)$ or
$(1-r_t(t+1))P_t(t+1)$ is $0$ because $P_t(-1)=P_t(t+1)=0$.
Outside $0\le k\le t+1$ both sides are $0$.
Under a fair start, $(X_i)$ and $(-X_i)$ have the same law, so
$P_t(k)=P_t(t-k)$. Taking second differences gives the same for $D_t$.

For \eqref{eq:D-rec}, take second differences in \eqref{eq:P-rec}.
The affine identity $r_t(k+1)-r_t(k)=a/t$ gives
\[
 2r_t(k+1)+r_t(k-2)=3r_t(k),\qquad
 r_t(k+1)+2r_t(k-2)=3r_t(k-1),
\]
so the coefficients of $P_t(k)$ and $P_t(k-1)$ agree with the
right side of~\eqref{eq:D-rec}. The outer two coefficients agree directly.
\end{proof}

The forward equation~\eqref{eq:P-rec} and a higher-difference
recurrence are also in Trinas and Valle~\cite[(2.2) and (4.10)]{tv2026}.
We use \eqref{eq:D-rec} to count sign changes. For $2\le k\le t-1$ its
coefficients are positive. At $k=1$ it reads
$D_{t+1}(1)=(1-r_t(2))D_t(1)+r_t(-1)D_t(0)$, where
$r_t(-1)=\varepsilon-a/t$ is negative for $t<a/\varepsilon$ and
$D_t(0)$ lies outside the list $D_t(1),\dots,D_t(t-1)$. So $s_{t+1}$ is
not a positive combination of neighbors. Hence, new changes are only
able to come in at the ends. Furthermore, Lemma~\ref{lem:interleave} shows that
positive sums of neighbors never add a change
(see~\cite{karlin1968,karlin1964} for general theory).

Let $z$ be a finite real list. Its sign pattern
$\operatorname{sp}(z)$ is an alternating list of signs generated by
removing the zeros, replacing each remaining entry by its sign, and
merging equal neighbors. For
two such lists we write $\sigma\preceq\tau$ if $\sigma$ is a subsequence of $\tau$. We put
$\Pi_5=(+,-,+,-,+)$ and $\Pi_3=(+,-,+)$.

\begin{lemma}[positive sums of neighbors]\label{lem:interleave}
Let $\zeta=(\zeta_1,\dots,\zeta_n)$ be real, and let
$\eta_i=A_i\zeta_i+B_i\zeta_{i+1}$ for $1\le i\le n-1$ with $A_i\ge0$ and
$B_i\ge0$. Let $w=(\zeta_1,\eta_1,\zeta_2,\eta_2,\dots,\eta_{n-1},\zeta_n)$.
Then $\operatorname{sp}(w)=\operatorname{sp}(\zeta)$ and
$\operatorname{sp}(\eta)\preceq\operatorname{sp}(\zeta)$.
\end{lemma}

\begin{proof}
Insert each $\eta_i$ between its two neighbors in $\zeta$. If $\eta_i=0$, the
sign pattern does not change. If $\eta_i>0$, then $\zeta_i>0$ or
$\zeta_{i+1}>0$, so $\eta_i$ is placed next to an entry of its own sign, and
the pattern does not change. The same holds for $\eta_i<0$, proving the first claim.
For the second, obtain $\eta$ from $w$ by deleting the entries
$\zeta_i$. Deletion can only remove signs from the pattern.
\end{proof}

We also use the proof in the following form. Each nonzero $\eta_i$
has the sign of $\zeta_i$ or of $\zeta_{i+1}$. We call such an entry
a parent of $\eta_i$. Two entries of $\eta$ with opposite signs have
different parents, so their parents come in the same order.

The recurrence does not control the values that arrive at the two
ends, so we carry them along from the start. Fix the final length $N$
and let $t<N$. To the left of the differences of time $t$ we write the
left end values of the times $N,N-1,\dots,t+1$, and to the right we
write the right end values of the times $t+1,\dots,N$. We call this
list the augmented row (Figure~\ref{fig:augmented-rows}, drawn for
first differences). The end value that arrives at time $t+1$ is
already next to its place in this list, so the recurrence controls
the whole list. Proposition~\ref{prop:all-bias-modes} treats the
whole row and Proposition~\ref{prop:dominant-half} the half
$2k+1\ge t$. Theorem~\ref{thm:five-blocks} uses the same row for
second differences, where both histories consist of the values $s_t$.

Let $P$ be a law on the integers. We call a maximal run of
consecutive integers on which $P$ is constant a \emph{modal plateau}
of $P$ if $P$ is strictly smaller on each side of it. We call a triple
$i<j<k$ a \emph{strict valley} of $P$ if $P(j)<\min\{P(i),P(k)\}$, and
we use the same words for a function on the line and points $x<y<z$.

\begin{figure}[ht]
\centering
\small
\[
\begin{array}{c@{\quad}c@{\quad}c}
\underbrace{\ell_N,\ldots,\ell_{t+1}}_{\text{left boundary history}}&
\underbrace{d_t(0),\ldots,d_t(t-1)}_{\text{active row}}&
\underbrace{u_{t+1},\ldots,u_N}_{\text{right boundary history}}\\[9pt]
&\downarrow\ \text{copy or sum adjacent entries}&\\[9pt]
\ell_N,\ldots,\ell_{t+2}&
\ell_{t+1},d_{t+1}(1),\ldots,d_{t+1}(t-1),u_{t+1}&
u_{t+2},\ldots,u_N
\end{array}
\]
\caption{One step of the augmented rows for a final length $N$, with
$d_t(k)=P_t(k+1)-P_t(k)$, $\ell_t=d_t(0)$, and $u_t=d_t(t-1)$. The
active row gains an entry and the histories on each side lose one.
Each new end value was already in a history, so the bound for the
previous row already counts its sign.}
\label{fig:augmented-rows}
\end{figure}

\begin{lemma}[sign changes next to an end]\label{lem:level-one}
For every $t_1\ge2$ we have $\operatorname{sp}(s_2,s_3,\dots,s_{t_1})\preceq\Pi_3$.
\end{lemma}

We prove this in Appendix~\ref{app:signs-level}.

\begin{theorem}[at most five blocks]\label{thm:five-blocks}
Let $N\ge3$ and $0<\varepsilon<\frac12$. Then
\[
 \operatorname{sp}\bigl(D_N(1),D_N(2),\dots,D_N(N-1)\bigr)\preceq(+,-,+,-,+).
\]
So the second difference of $P_N$ changes sign at most four times on
$\{1,\dots,N-1\}$.
\end{theorem}

\begin{proof}
For $3\le t\le N$, we retain the endpoint history in the list $\mathcal Z_t$
\[
 \mathcal Z_t=\bigl(s_N,s_{N-1},\dots,s_{t+1},\ s_t,D_t(2),D_t(3),\dots,D_t(t-2),s_t,\
 s_{t+1},\dots,s_{N-1},s_N\bigr).
\]
It has $2N-t-1$ entries. Its middle part is
$(D_t(1),D_t(2),\dots,D_t(t-1))$, because $D_t(t-1)=D_t(1)=s_t$ by symmetry.
The list $\mathcal Z_3$ is $(s_N,\dots,s_4,s_3,s_3,s_4,\dots,s_N)$, and $\mathcal Z_N$ is the
list in the statement.

To apply Lemma~\ref{lem:interleave}, we check how
$\mathcal Z_{t+1}$ is obtained from $\mathcal Z_t$ for $3\le t\le N-1$. The first $N-t$ entries of
$\mathcal Z_{t+1}$ are $s_N,\dots,s_{t+1}$, and they are copies of the first $N-t$
entries of $\mathcal Z_t$. So $A_i=1$ and $B_i=0$ there. The next $t-2$ entries are
$D_{t+1}(j)$ for $2\le j\le t-1$. By \eqref{eq:D-rec},
\[
 D_{t+1}(j)=r_t(j-2)\,D_t(j-1)+\bigl(1-r_t(j+1)\bigr)D_t(j),
\]
and $D_t(j-1)$ and $D_t(j)$ are neighbors in the middle part of $\mathcal Z_t$, in the
required order. The coefficients satisfy $r_t(j-2)\ge\varepsilon$ and
$1-r_t(j+1)\ge1-\varepsilon-a=\varepsilon$, since $2\le j\le t-1$. The last
$N-t$ entries of $\mathcal Z_{t+1}$ are $s_{t+1},\dots,s_N$, and each is a copy of the
entry of $\mathcal Z_t$ one place to its right. So $A_i=0$ and $B_i=1$ there.

By Lemma~\ref{lem:interleave},
$\operatorname{sp}(\mathcal Z_N)\preceq\operatorname{sp}(\mathcal Z_{N-1})\preceq\dots\preceq\operatorname{sp}(\mathcal Z_3)$.
By Lemma~\ref{lem:level-one} the sign pattern $p$ of $(s_3,\dots,s_N)$ is a
subsequence of $\Pi_3$, so it is one of the lists $()$, $(+)$, $(-)$, $(+,-)$,
$(-,+)$, and $(+,-,+)$. The sign pattern of $\mathcal Z_3$ is $p$ reversed, followed by
$p$, with equal neighbors merged. In the six cases this is $()$, $(+)$, $(-)$,
$(-,+,-)$, $(+,-,+)$, and $(+,-,+,-,+)$, and each is a subsequence of $\Pi_5$.
\end{proof}

\subsection{First differences for a biased start}\label{sec:biased-finite}

\begin{proposition}[two modes for every initial bias]\label{prop:all-bias-modes}
Let $0<\varepsilon<1$, let $q\in[0,1]$ be the probability that the first
step is positive, and let $P_n(k)$ be the law of the number of positive
steps after $n\ge1$ steps. Extend $P_n$ by zero outside $\{0,\ldots,n\}$.
The first differences $d_n(k)=P_n(k+1)-P_n(k)$ satisfy
\[
 \operatorname{sp}\bigl(d_n(-1),d_n(0),\ldots,d_n(n)\bigr)
 \preceq (+,-,+,-).
\]
In particular, $P_n$ has at most two modal plateaus, including endpoint
modes. Moreover, there are no integers $i_0<\cdots<i_4$ for which
\[
 P_n(i_0)>P_n(i_1)<P_n(i_2)>P_n(i_3)<P_n(i_4).
\]
\end{proposition}

\begin{proof}
We first bound the sign pattern of the left end values, and then
that of the first augmented row when $q\le1/2$. Then we show that one
step of the walk adds no sign. Reflection gives the case $q\ge1/2$.
Put $a=1-2\varepsilon$ and $r_t(k)=\varepsilon+ak/t$. The forward equation
is
\[
 P_{t+1}(k)=(1-r_t(k))P_t(k)+r_t(k-1)P_t(k-1),
 \qquad P_1=(1-q,q).
\]
Subtracting adjacent equations gives
\begin{equation}\label{eq:bias-first-diff}
 d_{t+1}(k)=r_t(k-1)d_t(k-1)+(1-r_t(k+1))d_t(k),
 \qquad 1\le k\le t-1.
\end{equation}
Both coefficients are at least $\min(\varepsilon,1-\varepsilon)>0$,
since $r_t(j)$ is a convex combination of $\varepsilon$ and
$1-\varepsilon$ for $0\le j\le t$. It remains to check the two ends.

\emph{The left end values.} Suppose first that $q<1$. Since
$P_t(0)=(1-q)(1-\varepsilon)^{t-1}>0$, define
\[
 \xi_t=\frac{P_t(1)}{P_t(0)},\qquad
 \beta=\frac{1-2\varepsilon}{1-\varepsilon}<1,\qquad
 \alpha=\frac{\varepsilon}{1-\varepsilon}>0.
\]
The forward equation at $k=0$ and $k=1$ gives
\[
 \xi_{t+1}=(1-\beta/t)\xi_t+\alpha,\qquad \xi_1=\frac q{1-q}.
\]
Thus $\delta_t=\xi_{t+1}-\xi_t$ obeys
\begin{equation}\label{eq:bias-edge-increments}
 \delta_{t+1}=\frac{t-\beta}{t+1}\delta_t+\frac\alpha{t+1}.
\end{equation}
Once an increment is nonnegative, every later increment is positive.
Thus $\xi_t$ first decreases and then increases (either part may be empty). For $\ell_t=d_t(0)=P_t(0)(\xi_t-1)$ this proves
\[
 \operatorname{sp}(\ell_1,\ldots,\ell_N)\preceq(+,-,+).
\]
If $q\le1/2$, then $\xi_1\le1$ and the sharper bound is
$\operatorname{sp}(\ell_1,\ldots,\ell_N)\preceq(-,+)$.
In particular, a positive $\ell_i$ cannot be followed by a negative
$\ell_j$. Zeros do not affect either conclusion.

\emph{The right end values.} Write $u_t=d_t(t-1)$. Reflection of the walk gives
$u_t(q)=-\ell_t(1-q)$ and
$d_t^{(1-q)}(t-1-k)=-d_t^{(q)}(k)$.
For $0<q\le1/2$, the right boundary history therefore has pattern
contained in $(-,+,-)$. If $q=0$ then $A_t\le t-1$, so $P_t(t)=0$ and
$u_t=-P_t(t-1)\le0$.

\emph{The first augmented row.} Fix $N$ and, for $1\le t\le N$, form the augmented row
\[
 Z_t=(\ell_N,\ldots,\ell_{t+1},
       d_t(0),\ldots,d_t(t-1),
       u_{t+1},\ldots,u_N).
\]
Let $q\le1/2$. The first row $Z_1$ is the left history read
backwards, followed by the right history, and the two share the entry
$\ell_1=u_1=2q-1$. The left history has pattern contained in $(-,+)$,
so read backwards its pattern is contained in $(+,-)$. The right
history has pattern contained in $(-,+,-)$, or in $(-)$ when $q=0$.
Hence the pattern of $Z_1$ is contained in $(+,-,+,-)$. This is where
the four signs of the proposition come from.

\emph{One step.} Every entry of $Z_{t+1}$ is a nonnegative combination of adjacent entries
of $Z_t$. The interior entries use \eqref{eq:bias-first-diff}. The two new
boundary entries and the remaining future edge values are copied,
with coefficients $(1,0)$ or $(0,1)$. Number the augmented entries from left to
right. Output entry $i$ uses input entries $i$ and $i+1$, with a zero
coefficient when copied. Hence Lemma~\ref{lem:interleave} gives
$\operatorname{sp}(Z_{t+1})\preceq\operatorname{sp}(Z_t)$, and so
$\operatorname{sp}(Z_N)\preceq(+,-,+,-)$. Reflection and sign
reversal preserve this pattern, proving the claim also for $q\ge1/2$,
including $q=1$.

\emph{Conclusion.} Extension by zero adds $d_N(-1)=P_N(0)\ge0$ on the left and $d_N(N)=-P_N(N)\le0$ on the right. The pattern $(+,-,+,-)$ starts with $+$ and ends with $-$, so the
longer list is still contained in it.
Every modal plateau contributes a positive-to-negative change, so there
are at most two. If the statement's last inequalities held, telescoping
on each of their four intervals would give an adjacent difference
of the same strict sign. These four differences would
have signs $-,+,-,+$, a contradiction.
\end{proof}

The last claim of Proposition~\ref{prop:all-bias-modes} is the form
that passes to the limit. If the limit had such strict inequalities at
five points, the count law would have them at nearby lattice points
for large $n$, and a count of plateaus has no such property. The next proposition keeps only the half of the row that the
start favors. With no strict valley on that half, the law has two modes
only if some first difference on the other half is negative. For each
pair of points the difference of the two values of the law is affine
in $q$, and this is the reason why the biases with two modes form an
interval. The argument is in the proof of
Theorem~\ref{thm:biased-limit-modes}, which also uses the inequality
$P_n(n-k)\le P_n(k)$ and the sign of the difference at $q=1$.

\begin{proposition}[the dominant half]\label{prop:dominant-half}
Under the hypotheses of Proposition~\ref{prop:all-bias-modes}, suppose
$q\ge1/2$. Then $P_n(n-k)\le P_n(k)$ for $2k\ge n$, and the right half has
no strict valley. That is,
\[
 P_n(j)\ge\min\{P_n(i),P_n(k)\}
 \qquad\text{if }i<j<k\text{ and }2i+1\ge n.
\]
Equivalently, its first differences contain no negative value followed
by a positive one. Reflection gives the corresponding left-half result
when $q\le1/2$.
\end{proposition}

\begin{proof}
We first show by induction that the law is
at least its reflection on the upper half. Then we run the induction
of Proposition~\ref{prop:all-bias-modes} on the right half of the
augmented row
$Z_t=(\ell_N,\ldots,\ell_{t+1},d_t(0),\ldots,d_t(t-1),u_{t+1},\ldots,u_N)$
of that proof, where the first part is used at one entry only.
Write the bias explicitly as $P_{n,q}$ and put
\[
 \Delta_n(k)=P_{n,q}(k)-P_{n,1-q}(k)
       =P_{n,q}(k)-P_{n,q}(n-k).
\]
The forward equation gives
$\Delta_{n+1}(k)=(1-r_n(k))\Delta_n(k)+r_n(k-1)\Delta_n(k-1)$.
At $n=1$, $\Delta_1(1)=2q-1\ge0$. Suppose that $\Delta_n(k)\ge0$
for $2k\ge n$. At a point strictly above the new center, both parents lie in
the preceding upper half and both coefficients are nonnegative.
At an integer center $\Delta_{n+1}(k)=0$ by reflection. Terms outside the
support vanish. This proves the induction. In particular, the
central first difference of an odd row is nonnegative.

We now keep only the half favored by the start, namely the entries
of $Z_t$ with spatial index $2k+1\ge t$ together with the right
history. At $t=1$ this list is the central entry $2q-1\ge0$ followed
by the future right end values $u_s(q)=-\ell_s(1-q)$. Since
$1-q\le1/2$, its pattern is contained in $(+,-)$ by
\eqref{eq:bias-edge-increments}, so it has no negative entry followed
by a positive one. Suppose that such a pair first appears at time
$t+1$. An entry of the right history is a copy of an entry of the
right history at time $t$, so we may take the negative entry in the
active row, at index $i$. Choose a parent of the same sign for each of
the two entries. By the remark after
Lemma~\ref{lem:interleave} the first parent comes before the second.
The only entry of the half at time $t+1$ that uses an entry outside the
half at time $t$ is the central difference of an odd row, and it is
nonnegative by the first part. So it is not the negative entry. Hence
$2i>t$, and both entries that the negative entry uses satisfy
$2k+1\ge t$.
The two parents then form a negative entry followed by a positive one
in the half at time $t$, which is a contradiction. A strict valley
would contain two such differences by telescoping, which proves the
claim.
\end{proof}

\section{The shape of the limit density}\label{sec:bulk-applications}

This section has four parts. Section~\ref{sec:modes} passes the finite
bounds to the limit for a fair start and states the threshold. The
next part does the same for every bias and finds the bias interval.
The third studies the bias boundary $b_a$, and the last the density at
the threshold. To pass a bound to the limit means the following. The
rescaled count law and its differences converge, locally uniformly, to
$g_{a,q}$ and its derivatives (Lemma~\ref{lem:derivative-limits}), and
a bound on the sign pattern of a sequence of functions holds for their
pointwise limit (Lemma~\ref{lem:sp-limit}).

Fix $\varepsilon\in(0,\frac14)$, so $a=1-2\varepsilon\in(\frac12,1)$.
Theorem~\ref{thm:five-blocks} holds for $\varepsilon<\frac12$ and
Proposition~\ref{prop:all-bias-modes} for $\varepsilon<1$. The limit
density exists only in the superdiffusive range $\varepsilon<\frac14$.
Conditional on $X_1=1$, $S_k/k^a$ converges almost surely to
$\mathcal L_a$~\cite[Theorem~3.7]{bercu2018}. The fair-start limit
$\mathcal L^\circ_a$ has the law of $\varsigma\mathcal L_a$, where $\varsigma$ is a fair sign independent of
$\mathcal L_a$~\cite[(1.3) and (1.6)]{glr2026}. Both limits have positive,
bounded, and smooth densities~\cite[Theorem~1.3]{glr2026}, the
densities $h_a$ and $g_a$ of the introduction.
We also use two regularity facts, proved in
Lemma~\ref{lem:analytic-density} of Appendix~\ref{app:regularity}. For
every bias $q$, the limit density extends to an entire function, so it
is not constant on any interval. For every $k\ge0$, the derivatives
$\partial_u^k h_a(u)$ and $\partial_a\partial_u^k h_a(u)$ exist and are
jointly continuous in $(a,u)$.

\subsection{Modes and inflection points}\label{sec:modes}

A probability density $h$ on $\mathbb R$ is \emph{unimodal} if for some $m$
it is nondecreasing on $(-\infty,m]$ and nonincreasing on $[m,\infty)$.
Gu\'erin, Laulin, Raschel, and Simon ask for a unique local maximum~\cite{glrs2025},
so we must also exclude a flat top. The density $g_a$ mixes the
unimodal density of $\mathcal L_a$ with its reflection. Such a
mixture can have more than one peak~\cite{dj1988}. We first take a fair start.
Lemma~\ref{lem:analytic-density} excludes constant intervals,
so modes are points. We consider other biases afterward.

\begin{theorem}[at most two modes]\label{thm:two-modes}
Let $a\in(\frac12,1)$. Exactly one of the following holds.
\begin{enumerate}[(a)]
\item The density $g_a$ is nondecreasing on $(-\infty,0]$ and nonincreasing on
$[0,\infty)$.
\item There is $\kappa>0$ such that $g_a$ is nondecreasing on $[0,\kappa]$ and
nonincreasing on $[\kappa,\infty)$, and $g_a(0)<g_a(\kappa)$.
\end{enumerate}
In case (b), $g_a$ is nondecreasing on $(-\infty,-\kappa]$ and nonincreasing on
$[-\kappa,0]$, and it is not unimodal. All these monotonicities are strict.
The unique local maximum is $0$ in case (a), and the only local maxima
are $\pm \kappa$ in case (b).
In case (a), there is $d>0$ such that $g_a$ is strictly concave on
$(-d,d)$ and strictly convex on each tail. In case (b), there are
$0<d_1<\kappa<d_2$ such that $g_a$ is strictly convex on $(-d_1,d_1)$ and
on each tail beyond $\pm d_2$, and strictly concave on
$(-d_2,-d_1)$ and $(d_1,d_2)$. Thus it has exactly two inflection points
in case (a) and four in case (b).
\end{theorem}

We call the statement on convexity in Theorem~\ref{thm:two-modes} the
curvature partition. The proof of Theorem~\ref{thm:bias-continuity}
uses it for the sign of $g_a'$ on each side of the mode.

A real function $v$ has \emph{sign pattern at most} $\Pi$ if
$\operatorname{sp}(v(y_1),\ldots,v(y_j))\preceq\Pi$ for every
$y_1<\cdots<y_j$.

\begin{proof}
We first show that $g_a''$ has sign pattern at most $\Pi_5$. Then we show
that $g_a$ has no strict valley on $[0,\infty)$, which gives the two
cases. Finally the $\Pi_5$ bound and the behavior of $g_a$ at zero and
at infinity give the sign pattern of $g_a''$.

Fix $a\in(1/2,1)$. Lemma~\ref{lem:derivative-limits} at order two
passes the five-block bound of Theorem~\ref{thm:five-blocks} to
$g_a''$. Indeed, the interpolant's values are positive sums of
adjacent second differences, and its lattice fills every fixed
compact interval as $N\to\infty$. The four values $D_N(0)$, $D_N(-1)$,
$D_N(N)$, and $D_N(N+1)$ near the ends sit at the points
$(2j-N)N^{-a}$ with $j\in\{-1,0,N,N+1\}$, since $\nabla^2P_N(k)=D_N(k-1)$.
These points have absolute value at least $N^{1-a}$. Since $a<1$, these points leave every compact
interval. Hence on a fixed compact interval the interpolant uses only
interior second differences for large $N$. The lattice spacing is
$2N^{-a}$, so for large $N$ any fixed points $y_1<\cdots<y_j$ lie in
different lattice cells. Lemmas~\ref{lem:interleave}
and~\ref{lem:sp-limit} therefore give sign pattern at most $\Pi_5$.

\emph{Monotone shape.} Proposition~\ref{prop:dominant-half} and
Lemma~\ref{lem:derivative-limits}, with $r=0$, show that $g_a$
has no strict valley on $[0,\infty)$. Such a valley would persist
at nearby lattice points under locally uniform convergence.

If $g_a$ is nonincreasing on $[0,\infty)$, case (a)
holds. Otherwise choose $y_1<y_2$ with $g_a(y_1)<g_a(y_2)$.
Since $g_a$ is integrable, there are $y_3<y_4$, both beyond $y_2$, with
$g_a(y_3)>g_a(y_4)$. The no-valley property excludes any later increase,
so $g_a$ attains a maximum on $[0,\infty)$ at some $\kappa$.
The absence of a strict valley makes it nondecreasing before
$\kappa$ and nonincreasing afterward. Since $g_a(y_1)<g_a(y_2)$, we have
$\kappa>0$. If $g_a(0)=g_a(\kappa)$,
it would be constant on $[0,\kappa]$, contrary to analyticity
and the positive increment. Thus case (b) holds.
Evenness gives the left-half shape and excludes unimodality in (b).

Equality at two points of a monotone interval would force a constant
interval, contrary to Lemma~\ref{lem:analytic-density}. Thus the
monotonicities are strict and there is no additional local maximum.

\emph{Curvature.} The Fourier bound in
Lemma~\ref{lem:analytic-density} gives $g_a(y),g_a'(y)\to0$ as
$y\to\infty$. Evenness and the $\Pi_5$ bound give at most three sign
blocks for $g_a''$ on $(0,\infty)$. Its last sign must be positive. Otherwise $g_a'$ would eventually
be nonincreasing with limit zero, making $g_a$ eventually nondecreasing.
This contradicts $g_a>0$ and $g_a\to0$.
Analyticity and strict monotonicity imply that the first nonconstant
term in the even Taylor series of $g_a$ at zero is negative in case (a)
and positive in case (b). The same signs therefore hold for $g_a''$
near zero. In case (b), $g_a''$ must also be negative somewhere, since
$g_a'$ is positive near zero and negative somewhere beyond $\kappa$.
The positive-half-line sign patterns are consequently $(-,+)$ in
case (a) and $(+,-,+)$ in case (b). Analyticity excludes an interval
of zero curvature and makes the resulting convexities strict. In
case (b), $g_a'(d_1)>0>g_a'(d_2)$, so $d_1<\kappa<d_2$.
\end{proof}

\begin{theorem}[the unique transition]\label{thm:threshold}
There is exactly one $a_*\in(1/2,1)$ with $g_{a_*}''(0)=0$.
It lies in $(2/3,69/100)$, and
\[
 g_a''(0)<0\quad(1/2<a<a_*),\qquad
 g_a''(0)>0\quad(a_*<a<1).
\]
Moreover, $\partial_a g_a''(0)|_{a=a_*}>0$.
The density $g_a$ is unimodal for $a\le a_*$ and bimodal for
$a>a_*$. In particular, the critical density is unimodal.
\end{theorem}

The proof is in Appendix~\ref{app:analytic-threshold}. The interval
$[2/3,69/100]$ comes from the proof. The weight in the formula~\eqref{eq:an-curvature} for the
curvature has a factor with the exponent $2-3a$, which changes sign at
$a=2/3$. Let $\mathcal R(a)$ be $(1+a)^2$ times the mean of the
integration variable under this weight, so that $g_a''(0)$ has the
sign of $\mathcal R(a)-1$. For $a\le2/3$ the factor lowers the mean
and the curvature is negative. For $a\ge2/3$ it raises the mean, and
the lower bound that this gives for $\mathcal R(a)$ exceeds $1$ from
$a=69/100$ on. Between the two values the proof shows that $\mathcal R$ is
increasing. Numerically $a_*\approx0.66881$
(Remark~\ref{rem:numerics}).

\subsection{Arbitrary initial bias}

Recall that $g_{a,q}(x)=q h_a(x)+(1-q)h_a(-x)$ is the limit density
for the initial bias $q$, and that $g_{a,1/2}=g_a$.

\begin{theorem}[the interval of bimodal biases]\label{thm:biased-limit-modes}
For every $a\in(1/2,1)$ and $q\in[0,1]$, the density $g_{a,q}$ has at
most two local maxima. There is $b_a\in(0,1/2]$ such that it has exactly
two local maxima if and only if
\[
                 b_a<q<1-b_a.
\]
Outside this open interval, including its endpoints, the density is
unimodal with a unique local maximum. The interval is empty exactly when
$a\le a_*$, in which case $b_a=1/2$. Positivity of $b_a$ for
$a>a_*$ is proved in Theorem~\ref{thm:bias-continuity}.
At finite time, for every $0<\varepsilon<1$, $q\in[0,1]$, and
$N\ge1$, the count law has at most two modal plateaus, including its
endpoints.
\end{theorem}

\begin{proof}
\emph{The mode bound.} Proposition~\ref{prop:all-bias-modes}
gives the finite assertion and excludes two successive strict valleys.
The locally uniform density convergence in
Lemma~\ref{lem:derivative-limits}, with $r=0$, passes these
inequalities to $g_{a,q}$.
Lemma~\ref{lem:analytic-density} excludes constant intervals
for every bias. Hence three local maxima $x<y<z$
would give
$x<u<y<v<z$ with
\[
 g_{a,q}(u)<\min\{g_{a,q}(x),g_{a,q}(y)\},\qquad
 g_{a,q}(v)<\min\{g_{a,q}(y),g_{a,q}(z)\}.
\]
For each pair, take a one-sided neighborhood of the smaller maximum
toward the other one. Since the density is not constant on an
interval, some point there has strictly smaller value. This gives
the two forbidden valleys.

\emph{The bias interval.} Fix $a$, and write $f_q=g_{a,q}$ and $h=h_a$.
Proposition~\ref{prop:dominant-half} and the same density convergence show that,
for $q\ge1/2$, $f_q$ has no strict valley on $[0,\infty)$ and
$f_q(-x)\le f_q(x)$ for $x\ge0$.
Thus it has a strict valley if and only if
$f_q(x)>f_q(y)$ for some $x<y\le0$. We call the inequality
$f_q(-x)\le f_q(x)$ for $x\ge0$ reflection dominance, and we call a
pair $x<y\le0$ with $f_q(x)>f_q(y)$ a descent. If $(x,y)$ is a descent
then $f_q(-x)\ge f_q(x)>f_q(y)$, so $(x,y,-x)$ is a strict valley. Conversely, let
$(x',y',z')$ be a strict valley. Then $x'<0$. If $y'\le0$, the pair
$(x',y')$ is a descent. If $y'>0$, the triple $(0,y',z')$ is not a
strict valley, so $f_q(0)\le f_q(y')<f_q(x')$, and the pair $(x',0)$
is a descent.

At $q=1$, conditional unimodality~\cite{glrs2025} and reflection
dominance place the mode in $[0,\infty)$. A negative mode would have
a reflected point of at least the same height, contrary to uniqueness.
Thus $h$ is nondecreasing on $(-\infty,0]$. For fixed $x<y\le0$, the
difference $f_q(x)-f_q(y)$ is affine in $q$ and nonpositive at $q=1$.
If positive at $q\ge1/2$, it remains positive for every $q'\in[1/2,q]$.
Thus, within $[1/2,1]$, a valley persists when the bias decreases.
These biases form a relatively open set because the strict
inequalities at three fixed points survive a small change of bias.
The set excludes $1$, and reflection gives an open interval
$(b_a,1-b_a)$, possibly empty.

A strict valley gives two local maxima. Integrability and positivity
let us choose an outer point on each side below the corresponding
high point. We then maximize on the two resulting compact
intervals. Conversely, two local
maxima give a strict valley because there is no flat interval.
The mode bound excludes a third maximum. A continuous integrable
density without a strict valley is unimodal, and absence of flat
intervals makes its maximum unique. Finally the bias interval is
nonempty exactly when it contains $1/2$, so
Theorem~\ref{thm:threshold} identifies the empty case.
\end{proof}

This proves the initial-bias conjecture
of~\cite[Remark~2.1(b)]{glrs2025}, with its half-width $q(a)=1/2-b_a$.
The derivative-limit argument used here is proved in Appendix~\ref{app:modes}.

\subsection{The bias boundary}

For $a>a_*$, Theorem~\ref{thm:threshold} puts $g_a$ in case (b) of
Theorem~\ref{thm:two-modes}. We write $\kappa_a$ for its positive
mode $\kappa$.

To find the endpoint, we compare the conditional density's slopes
at reflected points. Their ratio tells us how much weight the
reflected density needs to create a decrease.

\begin{proposition}[a variational equation for the bias boundary]
\label{prop:bias-boundary}
For $a\in(1/2,1)$ put $U_a=\{t>0:h_a'(t)>0\}$. With the infimum
over an empty set interpreted as $+\infty$,
\begin{equation}\label{eq:bias-boundary}
 b_a=\min\left\{\frac12,\quad
       \inf_{t\in U_a}\frac{h_a'(-t)}{h_a'(t)+h_a'(-t)}\right\}.
\end{equation}
If the infimum is attained
at $t>0$ with $h_a'(t)>0$ and $b_a<1/2$, then, with $q=1-b_a$,
\[
 g_{a,q}'(-t)=g_{a,q}''(-t)=0.
\]
For $a>a_*$ this is the point where the smaller peak and the valley
next to it merge (Theorem~\ref{thm:bias-continuity}).
\end{proposition}

\begin{proof}
The proof of Theorem~\ref{thm:biased-limit-modes} shows that
$h_a'(-t)\ge0$ for $t\ge0$. For $q\ge1/2$ bimodality is equivalent
to a decrease somewhere on the negative half-line, hence to
$g_{a,q}'(-t)<0$ for some $t>0$. Here
\[
 g_{a,q}'(-t)=q h_a'(-t)-(1-q)h_a'(t).
\]
This expression can be negative only if $t\in U_a$, in which case
it is negative exactly when
$1-q>h_a'(-t)/(h_a'(t)+h_a'(-t))$. The inequality is strict for each $t$, so the set of
$q\ge1/2$ with two modes is the set where $1-q$ exceeds the infimum.
Comparing with $q<1-b_a$ gives \eqref{eq:bias-boundary}. At an interior minimizer the
ratio has derivative zero. Differentiating its numerator and
denominator gives the two displayed equations.
\end{proof}

\begin{theorem}[regularity of the bias boundary]\label{thm:bias-continuity}
The function $a\mapsto b_a$ is continuous on $(1/2,1)$, is continuously
differentiable on $(a_*,1)$, and extends continuously to $b_1=0$.
It is strictly decreasing on some interval immediately above $a_*$.
For every $a>a_*$, we have $0<b_a<1/2$,
and the infimum in
\eqref{eq:bias-boundary} has a unique minimizer $t_a\in(0,\kappa_a)$.
This minimizer depends continuously on $a$ throughout $(a_*,1)$.
\end{theorem}

\begin{proof}
We first reduce the infimum to $(0,\kappa_a)$ and write it as the
maximum of a ratio $\Lambda_a$. We show that $\Lambda_a$ has no strict
valley, which gives a unique maximizer and $b_a>0$. Then we prove
continuity above $a_*$ and at $a_*$, and finally we compute the
derivative and its sign near $a_*$.

\emph{Reduction.} For $a>a_*$, the curvature partition in Theorem~\ref{thm:two-modes}
gives $g_a'(t)>0$ on $(0,\kappa_a)$ and $g_a'(t)<0$ on
$(\kappa_a,\infty)$. For $t\ge\kappa_a$ in $U_a$, the ratio in
\eqref{eq:bias-boundary} is therefore at least $1/2$.
Since $b_a<1/2$, it is enough to search $(0,\kappa_a)$
(we still have to show the infimum is attained).
Put
\[
 \Sigma_a(t)=h_a'(t)+h_a'(-t),\qquad
 \Lambda_a(t)=\frac{g_a'(t)}{\Sigma_a(t)}\quad(0<t<\kappa_a).
\]
Lemmas~\ref{lem:balanced-start} and~\ref{lem:positive-conditional-mode}
give $h_a'(-t)>0$ for every $t\ge0$.
We therefore have $\Sigma_a>0$, $0<\Lambda_a<1/2$, and
$1/2-b_a=\sup_{0<t<\kappa_a}\Lambda_a(t)$. Also $\Lambda_a(0+)=0$.

\emph{No strict valley.} We first rule out a strict valley in the ratio. We claim that
\begin{equation}\label{eq:bias-quasiconcave}
 \Lambda_a(v)\ge\min\{\Lambda_a(u),\Lambda_a(w)\}
 \qquad(0<u<v<w<\kappa_a).
\end{equation}
Otherwise, choose $\sigma\in(0,1/2)$ strictly between the middle value
and both outer values. With $q=1/2+\sigma$,
$g_{a,q}'(-t)=\Sigma_a(t)(\sigma-\Lambda_a(t))$ then has signs $-,+,-$
at $-w,-v,-u$, and is positive sufficiently close to zero on the left.
Positivity of the density
and its zero limits at infinity supply a positive derivative farther
left and a negative derivative farther right. There would be three
local maxima, contrary to Theorem~\ref{thm:biased-limit-modes}.

\emph{The maximizer.} At $\kappa_a$, we have $\Sigma_a(\kappa_a)=2h_a'(-\kappa_a)>0$.
Thus $\Lambda_a$ extends continuously to $[0,\kappa_a]$, with value
zero at both endpoints. Its positive maximum is attained in the
interior. At that point $h_a'(-t)>0$, so the maximum is strictly
below $1/2$. Hence $b_a>0$. We need attainment here (pointwise strict
inequality on an open interval would not suffice).
Equation~\eqref{eq:bias-quasiconcave} and analyticity make this
maximum unique, since two maximizers would force a constant interval.
They also make $\Lambda_a$ strictly increasing before its maximizer and
strictly decreasing after it.

\emph{Continuity above $a_*$.} Lemma~\ref{lem:analytic-density} gives joint continuity of all the
derivatives used here. The signs of $g_a'$ on either side of its
unique positive zero give continuity of $\kappa_a$. Fix $a_0>a_*$
and choose
$0<l<t_{a_0}<r<\kappa_{a_0}$. The strict inequalities
$\Lambda_a(l),\Lambda_a(r)<\Lambda_a(t_{a_0})$ persist near $a_0$.
For $v<l$, apply \eqref{eq:bias-quasiconcave} to
$v,l,t_{a_0}$. Since $\Lambda_a(l)<\Lambda_a(t_{a_0})$, it gives
$\Lambda_a(v)\le \Lambda_a(l)$. Applying it to $t_{a_0},r,v$ gives
$\Lambda_a(v)\le \Lambda_a(r)$ for $v>r$.
Thus every nearby global maximum lies in the fixed interval $[l,r]$.
The ratios converge uniformly on this fixed compact interval,
so their maxima converge. Uniqueness gives convergence of the
maximizers too. This proves continuity above $a_*$.

\emph{Continuity at $a_*$.} At $a_*$, the curvature partition in Theorem~\ref{thm:two-modes}(a)
makes $g_{a_*}'$ strictly decrease from zero and then strictly increase
to zero. Thus $g_{a_*}'(x)<0$ for every $x>0$.
Joint continuity gives $g_a'(x)<0$ for all $a>a_*$ sufficiently close
to $a_*$, so $\kappa_a<x$. Since $x$ was arbitrary, $\kappa_a\to0$.
On the shrinking interval $[0,\kappa_a]$, joint continuity now gives
\[
 \sup_{0\le t\le\kappa_a}|g_a'(t)|\longrightarrow0,\qquad
 \inf_{0\le t\le\kappa_a}\Sigma_a(t)\longrightarrow2h_{a_*}'(0)>0,
\]
where the strict inequality is Lemma~\ref{lem:positive-conditional-mode}.
Consequently $1/2-b_a=\max \Lambda_a\to0$. Below the threshold $b_a=1/2$,
so this proves continuity at $a_*$ without the critical expansion.
Corollary~\ref{cor:bias-boundary-rate} gives the limit at one.

\emph{The derivative.} For differentiability, the same three-point argument puts all nearby
maximizers in one compact interval. If $M(a)=1/2-b_a$, maximality gives,
for $s>0$,
\[
 \Lambda_{a+s}(t_a)-\Lambda_a(t_a)
 \le M(a+s)-M(a)
 \le \Lambda_{a+s}(t_{a+s})-\Lambda_a(t_{a+s}).
\]
Divide by $s$ and use the continuity of $t_a$ and
$\partial_a\Lambda_a$ from Lemma~\ref{lem:analytic-density}. The argument
for $s<0$ is the same. Thus
\begin{equation}\label{eq:bias-envelope-derivative}
 -b_a'=
 \frac{\partial_a h_a'(t_a)\,h_a'(-t_a)
       -h_a'(t_a)\,\partial_a h_a'(-t_a)}
      {\{h_a'(t_a)+h_a'(-t_a)\}^2}.
\end{equation}
The right side is continuous. No derivative of $t_a$ is needed.

\emph{The sign near $a_*$.} Put $C_*=\partial_a g_a''(0)|_{a=a_*}>0$ and
$H_*=h_{a_*}'(0)>0$, by Theorem~\ref{thm:threshold} and
Lemma~\ref{lem:positive-conditional-mode}.
Extend $\Lambda_a=g_a'/\Sigma_a$ to a fixed neighborhood of zero, where
$\Sigma_a>0$ for $a$ close to $a_*$. Evenness and
Lemma~\ref{lem:analytic-density} give, jointly as $(a,t)\to(a_*,0)$,
\[
 g_a'(t)=o(t),\qquad \partial_a g_a'(t)=(C_*+o(1))t,\qquad
 \Sigma_a(t)=2H_*+o(1),\qquad \partial_a\Sigma_a(t)=O(1).
\]
For the first two relations, integrate $g_a''$ and $\partial_a g_a''$
from zero to $t$, using $g_{a_*}''(0)=0$.
It follows that
\[
 \partial_a\Lambda_a(t)=\left(\frac{C_*}{2H_*}+o(1)\right)t>0
\]
for all sufficiently small $t>0$ and $a$ sufficiently close to $a_*$.
Since $0<t_a<\kappa_a\to0$, the envelope formula gives $b_a'<0$
on some interval immediately above the threshold.
\end{proof}

\begin{corollary}[a rate for the bias boundary near one]\label{cor:bias-boundary-rate}
Put
\[
 \rho_a=\frac{\Gamma(1+a)^2}{(2a-1)\Gamma(2a)}-1.
\]
If $\rho_a<1/12$, then
\[
 b_a\le\frac{4\rho_a}{1-4\rho_a}.
\]
In particular, as $a\uparrow1$,
$b_a\le8(1-a)+O((1-a)^2)$.
\end{corollary}

\begin{proof}
The conditional limit has mean $m_a=1/\Gamma(1+a)>0$ and variance
$\rho_am_a^2$, by the moment formulas of~\cite[Theorem~1.4]{glr2026}.
Consider the three adjacent intervals of length $m_a$ centered at
$-m_a$, $0$, and $m_a$. Chebyshev's inequality puts at least
$1-4\rho_a$ of the conditional law in the right interval and at most
$4\rho_a$ in the middle. For a start to the left, reflect these
intervals. After mixing, the left interval has mass at least
$(1-q)(1-4\rho_a)$ and the right at least $q(1-4\rho_a)$.
The middle still has mass at most $4\rho_a$.

For $4\rho_a/(1-4\rho_a)<q<1-4\rho_a/(1-4\rho_a)$, the middle
has less mass than either outer interval. All three intervals have
the same length. So, by continuity, we can choose a point in each
with the middle value strictly smaller than the other two.
This rules out unimodality. The two-mode bound in
Theorem~\ref{thm:biased-limit-modes} gives exactly two modes, and
therefore the stated bound on $b_a$. Finally, the logarithmic derivative at $a=1$
of $1+\rho_a$ is $2\psi(2)-2-2\psi(2)=-2$.
Thus $\rho_a=2(1-a)+O((1-a)^2)$, which gives the linear bound.
\end{proof}

Here $\rho_a$ is the variance of $\mathcal L_a$ divided by the square
of its mean. The corollary is explicit for every $a$ with
$\rho_a<1/12$ and uses only two moments. The linear bound is not sharp. Conditioning on a long run of positive
first steps puts most of the conditional law in a window of width of
order $1-a$, and this gives a quadratic rate.

\begin{proposition}[a quadratic rate near one]\label{prop:quadratic-rate}
As $a\uparrow1$,
\[
 \limsup_{a\uparrow1}\frac{b_a}{(1-a)^2}\le82.
\]
\end{proposition}

The proof is in Appendix~\ref{app:quadratic-rate}. It compares the mass
of $g_{a,q}$ on two adjacent intervals of equal length to the left of
zero.

\subsection{The density at the threshold}\label{sec:critical-density}

How far apart are the two modes just after the transition?
We compute this from derivatives at the threshold. As in the proof
of Theorem~\ref{thm:bias-continuity}, write
\[
 C_* = \left.\partial_a g_a''(0)\right|_{a=a_*},\qquad
 B_*=-g_{a_*}^{(4)}(0),\qquad H_*=h_{a_*}'(0).
\]
These constants are positive by Theorem~\ref{thm:threshold} and
Lemmas~\ref{lem:analytic-quartic} and~\ref{lem:positive-conditional-mode}.
The coefficient $Q_*$ is defined in Proposition~\ref{prop:bias-cusp}.

\begin{proposition}[the critical shape]\label{prop:critical-shape}
The density at the threshold satisfies
\[
 g_{a_*}(u)=g_{a_*}(0)-\frac{B_*}{24}u^4+O(u^6)
 \qquad(u\to0).
\]
Put $K_*=\sqrt{6C_*/B_*}$. As $a\downarrow a_*$, the positive mode
$\kappa_a$ of $g_a$ satisfies
\[
 \kappa_a=K_*\sqrt{a-a_*}+o\bigl(\sqrt{a-a_*}\bigr).
\]
\end{proposition}

\begin{proof}
Lemma~\ref{lem:analytic-density} gives joint continuity of all spatial
derivatives and their first memory derivatives. Taylor remainders
through order four are therefore uniform on a neighborhood of $a_*$.
At criticality, evenness, $g_{a_*}''(0)=0$, and the analytic sign
$g_{a_*}^{(4)}(0)=-B_*<0$ from Lemma~\ref{lem:analytic-quartic}
give the first assertion.

Put $\delta=a-a_*>0$. Uniform Taylor expansion gives
\[
 \frac{g_a'(u)}u=C_*\delta+o(\delta)
                  -\frac{B_*+o(1)}6u^2+O(u^4).
\]
For $u=k\sqrt\delta$ this has the sign of $C_*-B_*k^2/6$
when $k$ is fixed away from $\sqrt{6C_*/B_*}$. The unique positive
mode lies between any two such choices with opposite signs, by
Theorem~\ref{thm:two-modes}. Squeezing gives the mode expansion.
\end{proof}

\begin{proposition}[how the bias interval opens]\label{prop:bias-cusp}
With $C_*,B_*,H_*$ as defined above, put
\[
 Q_* = \frac{\sqrt{2}\,C_*^{3/2}}{3H_*\sqrt{B_*}}.
\]
As $a\downarrow a_*$,
\begin{equation}\label{eq:bias-cusp}
 \frac12-b_a=Q_*(a-a_*)^{3/2}+o\bigl((a-a_*)^{3/2}\bigr).
\end{equation}
The unique minimizer $t_a$ of
Theorem~\ref{thm:bias-continuity} satisfies
\[
 t_a=\sqrt{\frac{2C_*}{B_*}}\sqrt{a-a_*}
       +o\bigl(\sqrt{a-a_*}\bigr)
     =\frac{\kappa_a}{\sqrt3}+o\bigl(\sqrt{a-a_*}\bigr).
\]
At $q=1-b_a$, the stationary point is $-t_a$.
\end{proposition}

\begin{proof}
Lemma~\ref{lem:analytic-density} supplies the required joint regularity.
Theorem~\ref{thm:threshold} and
Lemma~\ref{lem:positive-conditional-mode} give $C_*>0$ and $H_*>0$.

Write $\delta=a-a_*$. By \eqref{eq:bias-boundary},
\[
 \frac12-b_a=\max\left\{0,\sup_{t\in U_a}
     \frac{g_a'(t)}{h_a'(t)+h_a'(-t)}\right\}.
\]
Theorem~\ref{thm:two-modes} makes the numerator nonpositive for $t\ge \kappa_a$.
Proposition~\ref{prop:critical-shape} gives
$\kappa_a=O(\sqrt\delta)$, and $H_*>0$ makes both conditional derivatives
positive on this shrinking interval. Uniform Taylor expansion gives,
for bounded $v\ge0$,
\[
 \frac{g_a'(v\sqrt\delta)}{h_a'(v\sqrt\delta)+h_a'(-v\sqrt\delta)}
 =\frac{\delta^{3/2}}{2H_*}
       \left(C_*v-\frac{B_*v^3}{6}\right)+o(\delta^{3/2}).
\]
The cubic is largest at the unique positive point
$v=\sqrt{2C_*/B_*}$. We maximize over $[0,\kappa_a/\sqrt\delta]$.
Its right endpoint tends to $\sqrt{6C_*/B_*}$, which is past that
point, and \eqref{eq:bias-cusp} follows. The ratio is zero at both
ends, so its positive maximum is attained in $(0,\kappa_a)$. Now the
convergence is uniform and the cubic has only one maximum. Hence
every scaled maximizer tends to $\sqrt{2C_*/B_*}$, which is the claim
for $t_a$.
\end{proof}

\section{Computer assistance}\label{sec:verify}

Every theorem, proposition, and lemma in this paper has a written proof
that uses no computer. A computer gives only the values in
Remark~\ref{rem:numerics} and the curves in Figure~\ref{fig:phase}.

\begin{remark}[numerical values]\label{rem:numerics}
Numerically, $a_*\approx0.66881$, $K_*\approx7.2241$, and
$Q_*\approx9.8670$.
\end{remark}
 The repository holds interval-arithmetic
certificates for $a_*$ and for $C_*$, $B_*$, $H_*$, $K_*$, and $Q_*$,
computed with Arb~\cite{johansson2017}, together with their records,
inputs, and acceptance tests. The repository also holds Lean files.
They verify all finite sign arguments and the analytic threshold proof
and assume the classical weak limit as a named hypothesis. The
coverage map lists the formal statements for each result. None of the
paper's proofs relies on the certificates or on the Lean files.

\section{Open problems}\label{sec:open}

\paragraph{The full bias boundary.}
Is $b_a$ strictly decreasing on $(a_*,1)$? By
Theorem~\ref{thm:bias-continuity} it is decreasing near $a_*$, and
Figure~\ref{fig:phase} suggests it is decreasing on the whole
interval. By \eqref{eq:bias-envelope-derivative}, $b_a'<0$ if and only
if the numerator there is positive at $t_a$. We cannot show this for
all $a$.

\paragraph{The rate near one.}
Proposition~\ref{prop:quadratic-rate} gives $b_a=O((1-a)^2)$. Is there
a matching lower bound $b_a\ge c(1-a)^2$ with $c>0$?

\paragraph{A finite-length threshold.}
For even $N$, does $D_N(N/2)$, as a function of $a$, have exactly one
zero in $(0,1)$? If so, does this zero converge to $a_*$?
Theorem~\ref{thm:threshold} is about the limit and does not answer
this.

\paragraph{The Drezner--Farnum model and higher dimensions.}
For $p\ne1/2$ the Drezner--Farnum rule is not $r_t(k)$, and the limit
is not a mixture of one law and its
reflection~\cite{dreznerfarnum1993}. Hence our argument does not
apply. In higher dimensions the fixed-point theory gives
densities~\cite{glr2026}, but the sign counting needs a linear order.

\medskip
\begin{acks}
I thank Peter Kagey, whose Problem~131 of the Open Problem Collection
\cite{kagey131} led to the crossing question, Aidan Deshong for
helpful conversations, and H\'elio Trinas for reading the paper.
\end{acks}

\section*{Code and data}
The repository
\url{https://github.com/apemm/elephant-bias-interval/tree/v1.4}
has the sources of this paper, the scripts that compute the numerical
values in Remark~\ref{rem:numerics} with their interval certificates
and records, and the Lean files. It is archived at Zenodo,
\href{https://doi.org/10.5281/zenodo.24095945}{doi:10.5281/zenodo.24095945}. The file
\nolinkurl{certificates/data/critical_constants_certificate.json}
holds the enclosures of the constants.

\appendix
\section{Derivative limits and the bias interval}\label{app:modes}

This appendix has two parts. Lemma~\ref{lem:derivative-limits} and
Lemma~\ref{lem:sp-limit} pass the finite bounds to the limit. The
theorems of Section~\ref{sec:bulk-applications} rest on the first, and
Theorem~\ref{thm:two-modes} also uses the second.
Lemma~\ref{lem:balanced-start} proves
$h_a'(t)>0$ for $t<0$ and Lemma~\ref{lem:positive-conditional-mode}
proves $h_a'(0)>0$. Together they give $b_a>0$ in
Theorem~\ref{thm:main-limit}.

\subsection{Derivative limits}\label{app:derivative-limits}

For $r=0,1$ the bounds and the local limit below are in Trinas and
Valle~\cite[Proposition~2.1 and Theorem~3.1]{tv2026}. We need every
$r$, and a fixed initial history, so we give the proof.

\begin{lemma}[derivative limits]\label{lem:derivative-limits}
Fix $a\in(1/2,1)$ and a finite initial history. Extend $P_n$ by zero
outside $\{0,\ldots,n\}$. Suppose the rescaled
count law converges weakly to a law with smooth density $f$.
For each fixed integer $r\ge0$, linearly interpolate
\[
 \frac{n^{(r+1)a}}{2^{r+1}}\nabla^rP_n(k)
 \quad\hbox{at}\quad (2k-n-r)n^{-a},
 \qquad \nabla P_n(k)=P_n(k)-P_n(k-1).
\]
These interpolants converge locally uniformly to $f^{(r)}$.
\end{lemma}

\begin{proof}
Apply $\nabla^r$ to \eqref{eq:P-rec}. This gives
\[
 \nabla^rP_{n+1}(k)
 =\left(\frac{1+a}{2}-\frac{ak}{n}\right)\nabla^rP_n(k)
 +\left(\frac{1-a}{2}+\frac{a(k-r-1)}n\right)\nabla^rP_n(k-1).
\]
For large $n$ and $1\le k\le n+r$, both coefficients are nonnegative
with sum $1-(r+1)a/n$. For $k=0$ and $k=n+r+1$ one term is zero and
the other has coefficient $(1+a)/2<1-(r+1)a/n$. Hence
$\sup_k|\nabla^rP_{n+1}(k)|\le(1-(r+1)a/n)\sup_k|\nabla^rP_n(k)|$,
and iterating,
\[
 \sup_k|\nabla^rP_n(k)|\le C_{a,r}n^{-(r+1)a}.
\]
Using this at orders $r$ and $r+1$, the interpolants are uniformly
bounded and uniformly Lipschitz. So every subsequence has a further
subsequence that converges locally uniformly. To find the limit, sum
by parts against a smooth function with compact support and use weak
convergence. What comes out is the $r$th distributional derivative of
the limit law. Note that the lattice spacing is $2n^{-a}$, and the
normalization above is exactly its $-(r+1)$ power. Since $f$ is
smooth, the distributional derivative is $f^{(r)}$. All subsequences
have this same limit, and the lemma follows.
\end{proof}

\subsection{Sign patterns}

This part has the sign lemma behind
Theorem~\ref{thm:two-modes}. Notation is as in
Sections~\ref{sec:convex} and~\ref{sec:modes}.
Deleting entries can only remove signs from a list's pattern. An alternating list $\sigma$ satisfies
$\sigma\preceq\tau$ for every alternating list $\tau$ with the same first sign and at
least as many entries.

\begin{lemma}[pointwise limits]\label{lem:sp-limit}
Let $v_n$ and $v$ be real functions on $\mathbb R$ with $v_n(y)\to v(y)$ for
every $y$, and let $\Pi$ be an alternating list. If every $v_n$ has sign pattern at
most $\Pi$, then so has $v$.
\end{lemma}

\begin{proof}
Let $y_1<\dots<y_j$, and let $I$ be the set of indices $i$ with $v(y_i)\ne0$.
For all large $n$, $v_n(y_i)$ has the sign of $v(y_i)$ for every $i\in I$.
Then $\operatorname{sp}(v(y_1),\dots,v(y_j))$ is the pattern of the list
$(v_n(y_i))_{i\in I}$, which is a subsequence of
$\operatorname{sp}(v_n(y_1),\dots,v_n(y_j))\preceq\Pi$.
\end{proof}

\subsection{The conditional density increases at zero}

Theorem~\ref{thm:bias-continuity} needs $h_a'>0$ on the whole negative
half-line. Unimodality of $h_a$~\cite{glrs2025} and reflection
dominance give only $h_a'\ge0$ there.

\begin{lemma}[a balanced start]\label{lem:balanced-start}
For $1/2<a<1$ the density conditional on a positive first step
satisfies
\[
 h_a'(t)>0\qquad(t<0).
\]
The proof uses the walk whose first two steps are $+1,-1$. For
$0\le a\le1$ its count law at each length $n\ge2$ is log-concave, with
positive support $\{1,\ldots,n-1\}$. For $1/2<a<1$ its limiting density
$f_a^{+-}$ is symmetric and log-concave, with $(f_a^{+-})'(t)<0$ for
$t>0$.
\end{lemma}

\begin{proof}
Given $X_1=1$, the second step is $+1$ with probability $(1+a)/2$. So
$h_a$ is a mixture of the limit densities $f_a^{++}$ and $f_a^{+-}$ of
the starts $+,+$ and $+,-$. We show that $f_a^{+-}$ is strictly
increasing on the negative half-line and that $f_a^{++}$ is
nondecreasing there. We prove the log-concavity of the count law for the
start $+,-$, identify its limit, pass the log-concavity to the limit,
and then treat the start $+,+$.
We run the log-concavity induction of~\cite[Proposition~2.2]{glrs2025} at
the points where $A,B\ge d$ below. Starting
with one step in each direction keeps the support inside
$\{1,\ldots,n-1\}$, so the two ends of the support need no check. The count law given $X_1=1$ is not
log-concave~\cite{tv2026}, so this induction cannot be run from one
fixed step.
Write one update as
$g_k=Af_k+Bf_{k-1}$, where
\[
 A=\frac{1+a}{2}-\frac{ak}{n},\qquad
 B=\frac{1-a}{2}+\frac{a(k-1)}n,\qquad d=\frac an.
\]
For $u=f_{k-2}$, $v=f_{k-1}$, $w=f_k$, and $z=f_{k+1}$, expansion gives
\begin{align*}
 g_k^2-g_{k-1}g_{k+1}
 &=(A^2-d^2)(w^2-vz)+(B^2-d^2)(v^2-uw)\\
 &\quad +(A-d)(B-d)(wv-uz)+d^2(w-v)^2.
\end{align*}
For $2\le k\le n-1$ we have $A,B\ge d\ge0$.
The preceding row is log-concave with interval support, so all
three gaps are nonnegative. At the two boundaries one of the
neighboring output masses is zero, and there is nothing to check. We
start from the point mass at one and induct. The recurrence also
shows the masses are strictly positive on the support.

Next we identify the limit. Step $t+1$ copies or reverses an earlier
step $X_i$. For $t\ge1$, join step $t+1$ to $X_i$ and give the
resulting edge the sign $X_{t+1}X_i$. Doing so creates a random
recursive tree on the steps, where the edge signs are independent both
of the tree and of each other. Each of these edge signs is $+1$ with
probability $(1+a)/2$. Removing the edge from step $2$ to step $1$
splits the steps into the part of step $1$ and the part of step $2$. A
new step joins the part that contains the step it copies or reverses.
Since that step is chosen uniformly, the new step joins a part with
probability equal to the fraction of the earlier steps that lie in
that part. Hence the sizes of the two
parts form a P\'olya urn started from one ball of each color, and the
share of the part of step $1$ converges almost surely to a uniform
variable $U$. Given the sequence of parts that the steps join, each
part is a random recursive tree in its own time, with its own
independent edge signs. So the sum of the steps in the part of step
$1$ is an elephant walk with a positive first step, observed at the
size of that part. The sum in the part of step $2$ is minus an
independent copy of such a walk. Both sizes tend to infinity. By the
almost sure convergence of $S_k/k^a$ given $X_1=1$~\cite[Theorem~3.7]{bercu2018},
\[
 \frac{S_n}{n^a}\longrightarrow
 \mathcal L_a^{+-}:=U^a\mathcal L_a'-(1-U)^a\mathcal L_a''
 \qquad\text{almost surely},
\]
where $\mathcal L_a'$ and $\mathcal L_a''$ are copies of $\mathcal L_a$
and $U$, $\mathcal L_a'$, $\mathcal L_a''$ are independent. With two
positive first steps the same split gives the limit
$\mathcal L_a^{++}=U^a\mathcal L_a'+(1-U)^a\mathcal L_a''$.

These limits have entire densities. Let $\phi_a$ be the characteristic
function of $\mathcal L_a$, and take $1<r<1/a$. Since $ra<1$, we have
$v^{ra}+(1-v)^{ra}\ge1$ for $0\le v\le1$. Hence the bound of
Lemma~\ref{lem:analytic-density} gives
\[
 \bigl|\E e^{it\mathcal L_a^{+\pm}}\bigr|
 \le\int_0^1|\phi_a(tv^a)|\,|\phi_a(t(1-v)^a)|\,dv
 \le4e^{-c|t|^r}.
\]
So, exactly as in the proof of that lemma, Fourier inversion gives
densities $f_a^{+-}$ and $f_a^{++}$ that extend to entire functions.
In particular they are smooth, so Lemma~\ref{lem:derivative-limits}
applies to both starts.

Now interpolate the rescaled lattice masses log-linearly. The result
is log-concave. By Lemma~\ref{lem:derivative-limits} with $r=0$ the
rescaled masses converge locally uniformly to $f_a^{+-}$. The
log-linear interpolant lies between adjacent masses, so it has the
same limit. A pointwise limit of log-concave functions is log-concave.
Hence $f_a^{+-}$ is log-concave. Symmetry follows from the
representation of $\mathcal L_a^{+-}$. The law of $\mathcal L_a^{+-}$
is a mixture over $U$ of convolutions of two positive densities, since
$h_a>0$. So $f_a^{+-}>0$. Now $(\log f_a^{+-})'$ is nonincreasing and is zero at zero.
Suppose it were also zero at some $t>0$. Then the density would be
constant on $[0,t]$, which analyticity rules out. Hence
$(f_a^{+-})'<0$ for $t>0$.

Next take two positive initial steps. We claim the first differences
are nonnegative on the negative half-line, and use induction from the
point mass at count two. Each entry on the negative half-line has
both parents in the lower half of the previous row. The one exception is the central
difference of an odd row, and there reflection dominance gives the
sign. The induction of Proposition~\ref{prop:dominant-half} runs in
the same way from the point mass at count two. The left end causes no
trouble, because the support starts at count two and the differences
at the left end are nonnegative. The derivative limit of Lemma~\ref{lem:derivative-limits} then
gives $(f_a^{++})'(t)\ge0$ for $t<0$.
Finally
$h_a=((1+a)/2)f_a^{++}+((1-a)/2)f_a^{+-}$.
The second summand has strictly positive derivative there.
\end{proof}

Gu\'erin, Laulin, Raschel, and Simon prove positivity of the conditional
mode for memory sufficiently close to one~\cite[Remark~2.1(a)]{glrs2025}.
We need it for all $a\in(1/2,1)$, in the stronger form of a positive
slope at zero.

\begin{lemma}[the conditional mode is positive]\label{lem:positive-conditional-mode}
For every $1/2<a<1$, we have $h_a'(0)>0$. In particular the unique
conditional mode lies strictly to the right of zero.
\end{lemma}

\begin{proof}
For $z=1$ or $z=2$ the number $N^{2a}e_N(z)/z$ defined below is a
rescaled central first difference, and it tends to $4h_a'(0)$. We show
that its lower limit is positive. After division by $z$ and a product
normalization, the reflection differences become averages of earlier
rows along a backward chain. The chain stays in a diffusive band with
probability above one half, and on the band the values at a fixed row
$n_0$ are bounded below. Write $f_n(k)=\Prob_+(A_n=k)$, with $\Prob_+$ the law given $X_1=1$, and,
for $z\equiv n\pmod2$,
$e_n(z)=f_n((n+z)/2)-f_n((n-z)/2)$.
With a positive start, the reflection induction of
Proposition~\ref{prop:dominant-half} is strict for $0<z\le n$. Indeed
$e_1(1)=1$, and each noncentral entry of the upper half has a
positive parent with a positive coefficient. The recurrence is
\begin{gather*}
 e_{n+1}(z)=\alpha_n(z)e_n(z+1)+\beta_n(z)e_n(z-1),\\
 \alpha_n(z)=\frac12-\frac{a(z+1)}{2n},\qquad
 \beta_n(z)=\frac12+\frac{a(z-1)}{2n}.
\end{gather*}
Since $\alpha_n(z)(z+1)+\beta_n(z)(z-1)=z(1-2a/n)$, the function $z$
is reproduced by the recurrence up to the factor $1-2a/n$. So dividing
by $z$ and by these factors turns the recurrence into an average, and
$e_N(z)/z$ is the slope of the count law at the center.
Fix $n_0>4a$. Divide by $z$ and multiply by
$\prod_{j=n_0}^{n-1}(1-2a/j)^{-1}$. Call the new values $v_n(z)$.
They satisfy a backward averaging equation, with weights
\[
 w_+=\frac{\alpha_n(z)(z+1)}{z(1-2a/n)},\qquad
 w_-=\frac{\beta_n(z)(z-1)}{z(1-2a/n)}.
\]
In the band we use below, these weights are nonnegative and add up to
one. At $z=1$ the parent $0$ gets weight zero. Let $z_n$ be the
backward chain with these weights. Where the weights are nonnegative,
$v_n(z_n)$ is a martingale in the backward time. Then
\[
 \E[z_n^2\mid z_{n+1}=z]
 =\frac{n-4a}{n-2a}z^2+\frac{3n-4a}{n-2a}.
\]
With $\gamma_n=\Gamma(n-2a)/\Gamma(n-4a)\sim n^{2a}$, the process
$Y_n=z_n^2/\gamma_n$ has backward drift
$\delta_n=(3n-4a)/\{(n-2a)\gamma_n\}$. Note that
$\sum_{n\ge n_0}\delta_n\to0$ as $n_0\to\infty$.

Fix $R>0$. We start the chain at $z_N=1$ if $N$ is odd and at $z_N=2$
if $N$ is even. We stop it at row $n_0$, or when it leaves the band
$0<z/\sqrt{\gamma_n}\le R$. Since $\sqrt{\gamma_n}\sim n^a$ and $a<1$,
we can take $n_0$ so large that $R\sqrt{\gamma_n}+1\le n$ for all
$n\ge n_0$. Then every state in the band and its two parents lie in
the support, and $\alpha_n(z),\beta_n(z)\ge0$ there. The chain uses the
weights only at band states before it stops, so we need their
nonnegativity only there. The chain moves by $\pm1$, and it can reach
$0$ only from $z=1$, where the parent $0$ has weight zero. So it
cannot leave the band through the bottom, and on leaving the band
$Y>R^2$. Up to the stopping time, the mean of $Y$ is at most
$z_N^2/\gamma_N+\sum_{n\ge n_0}\delta_n$. Taking $n_0$ and then $N$
large, Markov's inequality shows that the chain leaves the band with
probability below one half. Every stopped value of $v$ is positive,
by the strict reflection induction. The band states of row $n_0$ are
the finitely many $z$ with $0<z\le R\sqrt{\gamma_{n_0}}$, and at row
$n_0$ the normalizing product is empty, so $v_{n_0}(z)=e_{n_0}(z)/z$.
Their minimum $m_0>0$ depends only on $n_0$ and $R$, not on $N$. The
chain stops in this row with probability above one half. The
stopping time is bounded, so $v_N(z_N)$ is the mean of the stopped
value, and $v_N(z_N)\ge m_0/2$. The
product normalizer is a positive constant times $N^{2a}$ in the
limit. Hence
\[
 \liminf_N N^{2a}\frac{e_N(z_N)}{z_N}>0.
\]
By Lemma~\ref{lem:derivative-limits} for the conditional start,
the first differences converge, and this limit is $4h_a'(0)$. For
even $N$ we use that $e_N(2)$ is the sum of two adjacent central first
differences. So
$h_a'(0)>0$. Since $h_a$ is unimodal, its mode is to the right of
zero.
\end{proof}

\section{An analytic proof of the transition}\label{app:analytic-threshold}

Here we prove Theorem~\ref{thm:threshold}. The plan is to write the
sign of the central curvature as a comparison between positive
weights. We get the integral from the moment equation of Gu\'erin,
Laulin, Raschel, and Simon~\cite[(15)]{glrs2025}.

Let $E$ be an independent exponential variable of mean one. On the
imaginary axis the transform of $E^a\mathcal L_a$ has a closed
form~\cite[Proposition~3.2]{glrs2025}. The function $H_a$ defined next
parametrizes it in~\eqref{eq:an-param}, and its root $Y_a$ corresponds
to $t=\infty$. So we first collect what we need about $H_a$ and $Y_a$.

Put
\begin{equation}\label{eq:an-H}
 \chi_k=4^{-k}\binom{2k}{k},\qquad
 H_a(y)=1-\sum_{k\ge1}\frac{\chi_k}{2ak-1}y^k.
\end{equation}
Note that $H_a$ is strictly decreasing on $[0,1)$, and $H_1(y)=\sqrt{1-y}$.
For $a<1$ each coefficient $\chi_k/(2ak-1)$ is larger than its value at
$a=1$, so $H_a(1)<H_1(1)=0$.
(The series converges at $y=1$. Its coefficients are at most a
constant times $\chi_k/(2k-1)$, and these sum to one.) Hence $H_a$
has exactly one root $Y_a\in(0,1)$. Comparing coefficients, we also
get
\begin{equation}\label{eq:an-H-ode}
 H_a(y)-2ayH_a'(y)=(1-y)^{-1/2}.
\end{equation}
We have $\partial_aH_a(y)>0$ for $y>0$ and $H_a'(Y_a)<0$. By the
implicit function theorem, $Y_a$ is smooth and $Y_a'>0$.

\begin{lemma}[a weighted formula for the central curvature]\label{lem:curvature-weight}
For $1/2<a<1$,
\begin{equation}\label{eq:an-curvature}
 g_a''(0)=-\frac{1}{\pi\Gamma(3-3a)}
 \int_0^{Y_a}y^{1/2}H_a(y)^{2-3a}
             \bigl[1-(1+a)^2y\bigr]\,d y.
\end{equation}
The integral is absolutely convergent.
\end{lemma}

\begin{proof}
We first find a closed form for the transform of the limit on the
imaginary axis. Then we write its Mellin transform as an integral and
continue that integral to the point that gives the curvature.
We multiply the limit by an independent exponential power to cancel
the gamma factors in its moments.
Let $E$ be an independent exponential random variable of mean one
and set $Z_a=E^a\mathcal L_a$. With the normalization of~\cite{glrs2025},
\[
 m_n=\frac{\Gamma(1+an)}{n!}\E \mathcal L_a^n,
 \qquad m_0=m_1=1.
\]
The bound $0<m_n\le C_a^n$ follows from
\cite[Lemma~2.10 and Corollary~2.14]{glr2026}. Applied to
the even moments, it gives $\E\exp(\delta|Z_a|)<\infty$ for small
$\delta>0$. Hence $M_a(z)=\E e^{zZ_a}$ is analytic on a strip
about the imaginary axis. Compare its power series with the moment
equation. We get
\begin{equation}\label{eq:an-moment-ode}
 M_a(z)+azM_a'(z)
 =\frac{1+a}{2}M_a(z)^2
  +\frac{1-a}{2}M_a(z)M_a(-z)
\end{equation}
first near zero, and then on the whole strip by analytic
continuation. Write $M_a(it)=A(t)+iB(t)$. Then
\begin{equation}\label{eq:an-AB}
 atA'=A^2-A-aB^2,\qquad
 atB'=((1+a)A-1)B.
\end{equation}

For $0<y<Y_a$ define
\begin{equation}\label{eq:an-param}
 t=\sqrt y\,H_a(y)^{-a},\qquad
 A(t)=H_a(y)\sqrt{1-y},\qquad B(t)=H_a(y)\sqrt y.
\end{equation}
On the imaginary axis, this is the closed form of
\cite[Proposition~3.2]{glrs2025}, with their $G(it)=i\sqrt y$
and $(G(it)/(it))^{1/a}=H_a(y)$. We check the continuation directly
from~\eqref{eq:an-AB}. As $y$ runs over $(0,Y_a)$, $t$ increases from
zero to infinity, and by~\eqref{eq:an-H-ode}
\[
 \frac{1}{t}\frac{dt}{dy}
 =\frac{1}{2yH_a(y)\sqrt{1-y}}.
\]
So \eqref{eq:an-param} satisfies~\eqref{eq:an-AB}.
Next, $H_a(0)=1$, so the analytic map $y\mapsto yH_a(y)^{-2a}$ has
derivative one at zero. From $t^2=yH_a(y)^{-2a}$ we see that $y$ is
analytic in $t^2$ near zero. Hence $A$ is even and analytic, and
$B=tH_a(y)^{1+a}$ is odd and analytic, with $A(0)=1$ and $B'(0)=1$.
These values fix both power series, because the coefficient
recurrence has denominator $na-1$ for even $n\ge2$ and $na-a$ for
odd $n\ge3$, and both are positive. So near zero the parametrization
equals $M_a(it)$. By uniqueness for the ODE~\eqref{eq:an-AB} on
$t>0$, they agree for all $t>0$.

Write $f(t)=\operatorname{Re}\E e^{it\mathcal L_a}$. This is the
characteristic function of the fair limit, so
$g_a''(0)=-\frac1\pi\int_0^\infty t^2f(t)\,dt$, which is
$-\mathcal F(3)/\pi$ in the notation below. So we need $\mathcal F$ at
$s=3$. Exponential randomization gives
\[
 A(t)=\int_0^\infty e^{-v}f(tv^a)\,d v.
\]
We use the Fourier decay in Lemma~\ref{lem:analytic-density},
which in particular makes every polynomial multiple of $f$
integrable at infinity. For $0<\operatorname{Re}s<1/a$, absolute Fubini gives
\begin{equation}\label{eq:an-mellin}
 J(s):=\int_0^\infty t^{s-1}A(t)\,d t
 =\Gamma(1-as)\mathcal F(s),\qquad
 \mathcal F(s)=\int_0^\infty t^{s-1}f(t)\,d t.
\end{equation}
The function $\mathcal F$ is holomorphic on $\operatorname{Re}s>0$.
The integral defining $J$ diverges at $s=3$, the value needed for curvature.
Two integrations by parts give a convergent expression there.

As $t\to\infty$ the function $A(t)$ decays only like $t^{-1/a}$, since
$A(t)=y^{1/(2a)}\sqrt{1-y}\,t^{-1/a}$ by~\eqref{eq:an-param} and
$y\to Y_a$. So $J(s)$ diverges at that end once
$\operatorname{Re}s\ge1/a$. We change to the variable $u=t^{-1/a}$,
which moves that end to $u=0$, where the function $w$ below is smooth.
Each integration by parts then gains a power of $u$. Set $\theta=\arcsin\sqrt y$ and
\[
 u=F(\theta)=\sin^{-1/a}\theta\,H_a(\sin^2\theta),\qquad
 w(u)=a\cos\theta\sin^{1/a}\theta.
\]
Here $\theta$ decreases from $\arcsin\sqrt{Y_a}$ to zero as
$u$ increases from zero to infinity. From~\eqref{eq:an-H-ode},
\[
 F'(\theta)=-\frac1a\sin^{-1-1/a}\theta,\qquad t=u^{-a}.
\]
Changing variables in $J$ gives
$J(s)=\int_0^\infty u^{-as}w(u)\,d u$.
Since $F'(\arcsin\sqrt{Y_a})\ne0$, the function $w$ is smooth
at zero. At infinity, $\theta\sim u^{-a}$, so its definition
and the derivative formulas below give
$w(u)=O(u^{-1})$, $w'(u)=O(u^{-2})$, and $w''(u)=O(u^{-3})$.
For $0<\operatorname{Re}s<1/a$, the boundary terms therefore vanish twice,
and~\eqref{eq:an-mellin} becomes
\begin{equation}\label{eq:an-mellin-continued}
 \mathcal F(s)=\frac{1}{\Gamma(3-as)}
       \int_0^\infty u^{2-as}w''(u)\,d u.
\end{equation}
The right side is holomorphic on $0<\operatorname{Re}s<3/a$ by the bounds
at infinity and boundedness of $w''$ at zero.
Thus~\eqref{eq:an-mellin-continued} holds throughout this larger
strip by the identity theorem, and $s=3$ is permitted.

Differentiating through $u=F(\theta)$ gives
\[
 \begin{split}
 w'(u)&=-a\sin^{2/a}\theta
                  [1-(1+a)\sin^2\theta],\\
 w''(u)&=2a\sin^{3/a}\theta\cos\theta
                  [1-(1+a)^2\sin^2\theta].
 \end{split}
\]
Insert the second expression in~\eqref{eq:an-mellin-continued}
at $s=3$ and return to $y=\sin^2\theta$.
The powers of $y$ reduce to $y^{1/2}$, giving the integral in
\eqref{eq:an-curvature}. Fourier inversion supplies the factor
$g_a''(0)=-\mathcal F(3)/\pi$.
Finally, $H_a$ has a simple zero at $Y_a$, and $2-3a>-1$,
so the endpoint singularity is integrable.
\end{proof}

\begin{proof}[Proof of Theorem~\ref{thm:threshold}]
\emph{Signs outside a short interval.}

Let $\mu_a$ be the probability law on $(0,Y_a)$ with density proportional to
$y^{1/2}H_a(y)^{2-3a}$, and put
\begin{equation}\label{eq:an-score}
 \mathcal R(a)=(1+a)^2\int y\,d\mu_a(y).
\end{equation}
Lemma~\ref{lem:curvature-weight} says that $g_a''(0)$ has the sign
of $\mathcal R(a)-1$.

When $a\le2/3$, the extra weight $H_a^{2-3a}$ decreases in $y$.
By Chebyshev's association inequality, it lowers the mean of the density proportional to
$y^{1/2}$, giving $\int y\,d\mu_a(y)\le3Y_a/5$.
At $a=2/3$, the first three terms of~\eqref{eq:an-H} give
\[
 H_{2/3}(3/5)
 \le1-\frac9{10}-\frac{81}{1000}-\frac9{400}
 =-\frac7{2000}<0.
\]
Since $Y_a$ increases, $Y_a\le Y_{2/3}<3/5$ for $a\le2/3$.
Thus $\mathcal R(a)<(5/3)^2(3/5)^2=1$.

When $a\ge2/3$, the extra weight increases, giving a mean of at least
$3Y_a/5$. We need the bound
$Y_{2/3}>117/200$. To prove it, note that $\chi_k$ decreases.
At $y=117/200$, the terms with $k\ge4$ are at most
\[
 \sum_{k\ge4}\frac{\chi_ky^k}{4k/3-1}
 \le\frac{\chi_4y^4}{(13/3)(1-y)}.
\]
Together with the first three terms this gives the exact bound
\[
 H_{2/3}(117/200)\ge\frac{929771}{135987200}>0.
\]
For every $a\ge69/100$, therefore,
\[
 \mathcal R(a)>
 \left(\frac{169}{100}\right)^2\frac35\frac{117}{200}
 =\frac{10024911}{10000000}>1.
\]
Any zero must lie between $2/3$ and $69/100$.

\emph{A single crossing.}

We now work on $I=[2/3,69/100]$. To compare densities on a fixed
interval, rescale by $y=Y_au$, so $u$ always lies in $(0,1)$.
We separate the endpoint zero by putting
\begin{equation}\label{eq:an-weights}
 w_k(a)=\frac{\chi_kY_a^k}{2ak-1},\qquad
 K_a(u)=\frac{H_a(Y_au)}{1-u}
       =\sum_{k\ge1}w_k(a)(1+u+\cdots+u^{k-1}).
\end{equation}
Since $\sum w_k=1$, we have $K_a\ge1$ and $\partial_uK_a\ge0$.
Writing $s_a=3a-2\in[0,7/100]$, the rescaled measure has density
proportional to
\begin{equation}\label{eq:an-rescaled-weight}
 W_a(u)=u^{1/2}(1-u)^{-s_a}K_a(u)^{-s_a}.
\end{equation}
To make the mean increase, we show that the score
$\partial_a\log W_a(u)$ is strictly increasing in $u$ and use
the covariance identity for the derivative of an expectation.

In the derivative of the score in $u$, the factor $(1-u)^{-s_a}$
contributes $3/(1-u)\ge3$. The kernel $K_a$ contributes at most
$63/32$, and a term of size at most $45/4$ that is multiplied by
$s_a$. So the argument needs $s_a<11/120$, and $s_a\le7/100$ on $I$
suffices. This is why the interval $I$ is short.

First, $Y_a<2/3$ on $I$. The first three terms give
\[
 H_{69/100}(2/3)\le-\frac{9767}{7087608}<0,
\]
and $Y_a$ increases. So every series and derivative below converges
uniformly on $I\times[0,1]$. The falling-factorial moments are
\[
 M_1=\sum k w_k,\qquad M_2=\sum k(k-1)w_k,\qquad
 M_3=\sum k(k-1)(k-2)w_k.
\]
Let $G(y)=(1-y)^{-1/2}$. Evaluate~\eqref{eq:an-H-ode} and its first
two derivatives at $Y=Y_a$. This gives
\begin{equation}\label{eq:an-factorial-moments}
 \begin{split}
 2aM_1&=G(Y),\\
 2aM_2&=YG'(Y)-(2a-1)M_1,\\
 2a(M_3+2M_2)&=Y^2G''(Y)+M_2.
 \end{split}
\end{equation}
In particular,
\begin{equation}\label{eq:an-Ku}
 0\le\partial_uK_a(u)\le\frac{M_2}{2}
 \le\frac{YG'(Y)}{4a}
 \le\frac{3\sqrt3}{8}<\frac{21}{32}.
\end{equation}

For the memory derivatives of the kernel, put $d_k=2k/(2ak-1)$. Differentiating $\sum w_k=1$ gives
\begin{equation}\label{eq:an-root-derivative}
 \frac{Y_a'}{Y_a}M_1=\sum d_kw_k\le6,
 \qquad \frac{Y_a'}{Y_a}\le\frac{12a}{G(Y)}.
\end{equation}
Here $d_k$ decreases in $k$ and $d_1\le6$ on $I$.
The logarithmic derivative $\partial_a\log w_k=kY_a'/Y_a-d_k$
strictly increases in $k$. Since $1+u+\cdots+u^{k-1}$
is nondecreasing in $k$, the covariance identity for
derivatives of expectations gives $\partial_aK_a(u)\ge0$.
The probability law here is formed by the $w_k$.

In $\partial_a\partial_uK_a$, discard the negative term
$-d_kw_k$ from $\partial_aw_k$, then bound the derivative
of $1+u+\cdots+u^{k-1}$ by $k(k-1)/2$. Equations
\eqref{eq:an-factorial-moments} and~\eqref{eq:an-root-derivative} give
\begin{align}
 \partial_a\partial_uK_a
 &\le\frac{Y_a'}{Y_a}\frac{M_3+2M_2}{2}\notag\\
 &\le3\left(\frac{M_2}{G(Y)}+
                       \frac{Y^2G''(Y)}{G(Y)}\right)
 \le3\left(\frac34+3\right)=\frac{45}{4}.
 \label{eq:an-Kau}
\end{align}
The last bounds use
$M_2/G(Y)\le Y/[4a(1-Y)]$ and
$Y^2G''(Y)/G(Y)=3Y^2/[4(1-Y)^2]$.
Since $K_a\ge1$, $\partial_aK_a\ge0$ and
$\partial_uK_a\ge0$, we obtain
\[
 \partial_a\left(\frac{\partial_uK_a}{K_a}\right)
 =\frac{\partial_a\partial_uK_a}{K_a}
  -\frac{\partial_uK_a\,\partial_aK_a}{K_a^2}
 \le\frac{45}{4}.
\]
In~\eqref{eq:an-rescaled-weight}, the endpoint term dominates the kernel
terms, giving the strict bound
\begin{align}
 \partial_u\partial_a\log W_a(u)
 &=\frac3{1-u}-3\frac{\partial_uK_a}{K_a}
    -s_a\partial_a\left(\frac{\partial_uK_a}{K_a}\right)\notag\\
 &\ge3-\frac{63}{32}-\frac7{100}\frac{45}{4}
   =\frac{39}{160}>0.
 \label{eq:an-increasing-score}
\end{align}
The score $\partial_a\log W_a$ is therefore strictly increasing
in $u$. If $\nu_a$ denotes the normalized measure in
\eqref{eq:an-rescaled-weight}, differentiation gives
\[
 \frac{d}{da}\int u\,d\nu_a(u)
 =\operatorname{Cov}_{\nu_a}\bigl(u,\partial_a\log W_a(u)\bigr)>0.
\]
This is the strict form of Chebyshev's association inequality. We may differentiate under the
integral, because on this compact parameter interval the integrand
is dominated by $C(1-u)^{-7/100}(1+|\log(1-u)|)$, which is integrable.

It follows that
\[
 \mathcal R(a)=(1+a)^2Y_a\int u\,d\nu_a(u)
\]
has strictly positive derivative in the interior of $I$. By the signs
at the ends of $I$ and the intermediate value theorem, $\mathcal R-1$
has exactly one zero in $I$. We already know it has none outside $I$.
Now rewrite~\eqref{eq:an-curvature} as
\[
 g_a''(0)=\frac{\mathcal Z_a}{\pi\Gamma(3-3a)}\,[\mathcal R(a)-1],
 \qquad \mathcal Z_a=\int_0^{Y_a}y^{1/2}H_a(y)^{2-3a}\,d y>0.
\]
The same rescaling and domination show that the prefactor is
continuously differentiable near the zero. At the zero
$\mathcal R-1$ vanishes, so only $\mathcal R'$ contributes, and
$\partial_ag_a''(0)=\mathcal Z_a\mathcal R'(a)/[\pi\Gamma(3-3a)]>0$.
By Theorem~\ref{thm:two-modes}, $g_a$ is then unimodal for $a<a_*$
and bimodal for $a>a_*$. For $a=a_*$ itself, let $a\uparrow a_*$. The
inequality $g_a(x)\ge g_a(y)$ for $0\le x<y$ survives in the limit by
joint continuity. So $g_{a_*}$ is nonincreasing on $[0,\infty)$, and
since it is not constant on any interval, its only mode is zero.
\end{proof}

With the same bounds on $Y_a$ we can also get the sign of the
fourth derivative.

\begin{lemma}[the quartic sign]\label{lem:analytic-quartic}
For $2/3\le a\le69/100$, we have $g_a^{(4)}(0)<0$.
In particular, $B_*=-g_{a_*}^{(4)}(0)>0$.
\end{lemma}

\begin{proof}
We integrate by parts once more in the Mellin calculation.
With $u,w$ as in the proof of Lemma~\ref{lem:curvature-weight}, put
\[
 c_a=(1+a)^2,\qquad
 P_a(y)=\frac3a(1-y)(1-c_ay)-y(1-c_ay)-2c_ay(1-y).
\]
Since $dy/du=-2ay^{1+1/(2a)}\sqrt{1-y}$, differentiation gives
$w'''(u)=-2a^2y^{2/a}P_a(y)$.
This derivative is bounded at zero and is $O(u^{-4})$ at infinity.
The boundary terms vanish for $0<\operatorname{Re}s<3/a$, so
\eqref{eq:an-mellin-continued} gives
\[
 \mathcal F(s)=-\frac1{\Gamma(4-as)}
                \int_0^\infty u^{3-as}w'''(u)\,d u.
\]
The right side is holomorphic for $0<\operatorname{Re}s<4/a$.
For $1/2<a<4/5$, we may therefore set $s=5$ and use
$g_a^{(4)}(0)=\mathcal F(5)/\pi$ to obtain
\begin{equation}\label{eq:an-fourth-derivative}
 g_a^{(4)}(0)=\frac{a}{\pi\Gamma(4-5a)}
 \int_0^{Y_a}y^{3/2}(1-y)^{-1/2}H_a(y)^{3-5a}P_a(y)\,d y.
\end{equation}
The integral converges absolutely because $H_a$ has a simple zero
and $3-5a>-1$.

Now let $a\in[2/3,69/100]$. The root bounds in the preceding proof
give $1/2<Y_a<2/3$. For the coefficient comparison only, we use
the larger interval $2/3\le a\le7/10$. Expanding $P_a$ gives, for $y\ge0$,
\[
 \begin{split}
 P_a(y)&=\frac3a-\left(2a^2+7a+9+\frac6a\right)y
                    +\left(3a^2+9a+9+\frac3a\right)y^2\\
       &\le\frac92-23y+22y^2=:\widehat P(y).
 \end{split}
\]
Here termwise lower bounds give
$8/9+14/3+9+60/7=1457/63>23$ for the coefficient of $-y$,
and termwise upper bounds give
$147/100+63/10+9+9/2=2127/100<22$ for that of $y^2$.
Since $\widehat P(0)>0$, $\widehat P(1/2)=-3/2$, and $\widehat P(2/3)=-19/18$,
the quadratic $\widehat P$ has exactly one root $y_0$ in $(0,Y_a)$,
with $\widehat P>0$ before it and $\widehat P<0$ after it. Integrating
directly, for $1/2\le Y\le2/3$,
\[
 Y^{-5/2}\int_0^Y y^{3/2}\widehat P(y)\,d y
 =\frac95-\frac{46}{7}Y+\frac{44}{9}Y^2
 \le-\frac{83}{315}<0.
\]
The middle expression is convex in $Y$, so it suffices to check
$Y=1/2$ and $Y=2/3$.

The weight $\omega_a(y)=(1-y)^{-1/2}H_a(y)^{3-5a}$ is what remains,
and it is positive and increasing on $(0,Y_a)$. Since $\widehat P$
swaps signs only at $y_0$, we have
$\omega_a(y)\widehat P(y)\le \omega_a(y_0)\widehat P(y)$ on the whole
interval. Hence
\[
 \int_0^{Y_a}y^{3/2}\omega_a(y)P_a(y)\,d y
 \le \omega_a(y_0)\int_0^{Y_a}y^{3/2}\widehat P(y)\,d y<0.
\]
On $[2/3,69/100]$ the prefactor in~\eqref{eq:an-fourth-derivative} is
positive, and the lemma follows.
\end{proof}

\section{Signs at the first lattice point}\label{app:signs}

\label{app:signs-level}
We prove Lemma~\ref{lem:level-one}. The value $s_t$ satisfies a
linear recursion in $t$ with a forcing term. After division by a
positive product, $s_t$ is a partial sum of $s_2$ and the forcing
terms. The forcing terms change sign at most twice, so the partial
sums change sign at most twice. Put $\xi_t=P_t(1)/P_t(0)$, the ratio
of the proof of Proposition~\ref{prop:all-bias-modes} for a fair
start. By \eqref{eq:P-rec}, $P_{t+1}(0)=(1-\varepsilon)P_t(0)$
and $P_{t+1}(1)=(1-\varepsilon-\frac at)P_t(1)+\varepsilon P_t(0)$. So
\[
 \xi_{t+1}-\xi_t=\frac1{1-\varepsilon}\Bigl(\varepsilon-\frac{a\xi_t}t\Bigr),
 \qquad\xi_2=\frac{2\varepsilon}{1-\varepsilon}.
\]
We claim that $\xi_t<\varepsilon t/a$ for $t\ge2$. For $t=2$ this holds
because $a<1-\varepsilon$. If it holds for $t$, then
$0<\xi_{t+1}-\xi_t\le\frac\varepsilon{1-\varepsilon}<\frac\varepsilon a$,
and it holds for $t+1$. So $\xi_t$ is strictly increasing for $t\ge2$, and
$\xi_t<1$ whenever $t\le a/\varepsilon$.

Now \eqref{eq:D-rec} with $k=1$, and $D_t(0)=P_t(1)-2P_t(0)$, give for $t\ge2$
\[
 s_{t+1}=\bigl(1-r_t(2)\bigr)s_t+\sigma_t,\qquad
 \sigma_t=\Bigl(\varepsilon-\frac at\Bigr)(\xi_t-2)\,P_t(0).
\]
Here $1-r_t(2)\ge1-\varepsilon-a=\varepsilon>0$. If $t<a/\varepsilon$, both
factors of $\sigma_t$ are negative and $\sigma_t>0$. If $t>a/\varepsilon$, the
first factor is positive, and the second is first negative and then positive,
because $\xi_t$ increases. So
$\operatorname{sp}(\sigma_2,\sigma_3,\dots)\preceq\Pi_3$. To remove the positive multiplier, let
$\Phi_t=\prod_{u=2}^{t-1}(1-r_u(2))>0$. Then
\[
 \frac{s_t}{\Phi_t}=s_2+\sum_{u=2}^{t-1}\frac{\sigma_u}{\Phi_{u+1}},\qquad
 s_2=1-3\varepsilon .
\]
If $\varepsilon<\frac13$, the terms $s_2,\sigma_2,\sigma_3,\dots$ are first
nonnegative, then nonpositive, then nonnegative. Their partial sums are
nonnegative and nondecreasing in the first stretch, nonincreasing in the
second, and nondecreasing in the third. So the sign pattern of the partial
sums is a subsequence of $\Pi_3$. If $\varepsilon\ge\frac13$, then
$a/\varepsilon\le1$, so the terms are first nonpositive and then nonnegative,
and the sign pattern of the partial sums is a subsequence of $(-,+)$. This
proves the claim. \qed

\section{Regularity of the limit densities}\label{app:regularity}

We prove the regularity used in Sections~\ref{sec:bulk-applications}
and~\ref{app:modes}.

\begin{lemma}[entire densities and parameter continuity]\label{lem:analytic-density}
For every $a\in(\frac12,1)$ and $1<r<1/a$, the characteristic function
$\phi_a$ of $\mathcal L_a$ satisfies
\[
 |\phi_a(t)|\le 2e^{-c|t|^r}\qquad(t\in\mathbb R)
\]
for some $c=c(a,r)>0$. The same $c$ works on a sufficiently small
neighborhood of any fixed $a$. For every compact
$J\subset(1/2,1)$ and every $\beta>0$, there is $C_{J,\beta}<\infty$ such that
\[
 |\partial_t\phi_a(t)|+|\partial_a\phi_a(t)|
 \le C_{J,\beta}(1+|t|)^{-\beta}
 \qquad(a\in J,\ t\in\mathbb R).
\]
For every integer $k\ge0$, the derivatives $\partial_u^k h_a(u)$ and
$\partial_a\partial_u^k h_a(u)$ exist and are jointly continuous in $(a,u)$.
For every $q\in[0,1]$, the density of the limit with first-step bias $q$
extends to an entire function and is not constant on an interval of
positive length.
\end{lemma}

\begin{proof}
After weak continuity in $a$, a fixed-point inequality gives an
exponential bound for the characteristic function, so the densities
are entire. Next the moment series gives analyticity in $a$
at bounded frequencies, and finally the same fixed-point argument
gives polynomial decay of the two derivatives.
We first check weak continuity in $a$ to obtain a Fourier bound local in memory. On a compact
$J\subset(1/2,1)$, the martingale
$\mathcal M_n(a)=\Gamma(n)S_n/\Gamma(n+a)$ of~\cite[(2.4), Appendix~A]{bercu2018} satisfies
\[
 \sup_{a\in J}\E|\mathcal L_a-\mathcal M_n(a)|^2
 \le C_J\sum_{k\ge n}k^{-2\min J}\longrightarrow0.
\]
The walk's centered next step has conditional variance at most one,
and the gamma ratios are uniformly bounded on $J$. The finite laws are
continuous in $a$, so $\phi_a(t)$ is jointly continuous, uniformly on
compact frequency intervals.

The fixed-point equation of~\cite[Theorem~1.1 and (2.10)]{glr2026} gives
\[
 |\phi_a(t)|\le\int_0^1|\phi_a(tv^a)|\,
                  |\phi_a(t(1-v)^a)|\,dv.
\]
Here $\phi_a$ is the characteristic function conditional on a positive
first step. The random sign disappears after taking absolute values, since
$|\phi_a(-s)|=|\phi_a(s)|$.
Fix $a_0$ and $1<r<1/a_0$, and choose $J$ about $a_0$ with $r\max J<1$.
Put $D(v)=v^{r\max J}+(1-v)^{r\max J}-1$. We have $D(v)>0$ for
$0<v<1$. By dominated convergence there is $K>1$ such that
\[
 R:=2\int_0^1e^{-K^r(\log2)D(v)}\,dv<1.
\]
The law of $\mathcal L_{a_0}$ has a density~\cite[Theorem~1.3]{glr2026}, so
the Riemann--Lebesgue lemma gives $T>0$ with
$|\phi_{a_0}(t)|<2^{-K^r}$ for $|t|\ge T$. Joint continuity, after
shrinking $J$, gives $|\phi_a(t)|\le2^{1-K^r}$ for $a\in J$ and
$T\le|t|\le KT$. Set
$c=(\log2)/T^r$ and $b(t)=2e^{-ct^r}$ for $t\ge0$.
On $[0,T]$ we have $b\ge1$, and on $[T,KT]$ we have
$b\ge2^{1-K^r}$. Thus $|\phi_a|\le b$ on $[0,KT]$.

For $t\ge KT$,
\[
 \frac{\int_0^1b(tv^a)b(t(1-v)^a)\,dv}{b(t)}
 \le2\int_0^1e^{-ct^rD(v)}\,dv\le R<1.
\]
Suppose the bound $|\phi_a|\le b$ stops at a finite
$t_*\ge KT$. By continuity it still holds on $[0,t_*]$.
But at $t_*$ the fixed-point inequality gives
$|\phi_a(t_*)|\le Rb(t_*)<b(t_*)$, which is strict. By continuity the
bound holds a little past $t_*$, a contradiction. Hence it holds on
$[0,\infty)$, and on $\mathbb R$ by evenness.

Since $r>1$, the Fourier integral converges absolutely and locally
uniformly for complex arguments. Indeed, for every $B<\infty$,
$\int_{\mathbb R}e^{B|t|}|\phi_a(t)|\,dt<\infty$.
It therefore defines an entire function.
After any fixed number of spatial derivatives, the bound on $J$
still supplies an integrable majorant, giving joint continuity by
dominated convergence. For initial bias $q$, the triangle inequality
applied to $q\phi_a(t)+(1-q)\phi_a(-t)$ gives the same bound.
A density of integral one cannot be constant on $\mathbb R$.
If it were constant on an interval, the identity theorem would
make it constant everywhere.

To differentiate in memory, we first work on bounded frequency intervals.
The moment recursion of~\cite[Theorem~1.4]{glr2026} is
\[
 m_0=m_1=1,\qquad
 m_n(z)=\frac{\sum_{j=1}^{n-1}\theta_j(z)m_j(z)m_{n-j}(z)}
                   {nz-\theta_n(z)},\qquad
 \theta_j(z)=\begin{cases}1&j\text{ even},\\z&j\text{ odd}.
                  \end{cases}
\]
Choose a closed complex disc about $a_0$ on which
$\operatorname{Re}z\ge\alpha>1/2$, $|z|<1$ and
$|\operatorname{Im}z|\le\rho$. All denominators for $n\ge2$ have
real part at least $n\alpha-1>0$.
Put $A=2\alpha-1$ and $B=A^{-1}$. Induction gives
\[
 |m_n(z)|\le AB^n\quad(n\ge1),\qquad
 \frac{(n-1)A^2}{n\alpha-1}\le A\quad(n\ge2).
\]
For $x>0$ and real $y$, the gamma product gives
\[
 \left|\frac{\Gamma(x)}{\Gamma(x+iy)}\right|^2
 =\prod_{j\ge0}\left(1+\frac{y^2}{(x+j)^2}\right)
 \le\exp\left(\frac{y^2}{x^2}+\frac{y^2}{x}\right).
\]
Taking $x=1+n\operatorname{Re}z$ and $y=n\operatorname{Im}z$ gives
\[
 \left|\frac{m_n(z)}{\Gamma(1+nz)}\right|
 \le\frac{C B_0^n}{\Gamma(1+n\alpha)}.
\]
Here $C,B_0$ are independent of $n,z$. We used monotonicity of
$\Gamma$ on $[2,\infty)$ for $n\ge2$, and enlarged $C$ for $n=0,1$.
Thus $\sum_{n\ge0}m_n(z)(it)^n/\Gamma(1+nz)$ converges locally
uniformly in complex $t$ and $z$ in the disc's interior.
For real $z=a$ this is $\phi_a(t)$, by the moment formula and
exponential moments of~\cite[Corollary~2.11]{glr2026}.
Cauchy's formula gives local joint continuity of every frequency
derivative and its memory derivative.

We now obtain bounds at large frequencies.
Shrink $J$ about $a_0$ so that it lies inside this disc and retains the
uniform exponential bound. For $t\ge0$ write
$f_a(t)=|\phi_a(t)|$, $u_a(t)=|\partial_t\phi_a(t)|$ and
$w_a(t)=|\partial_a\phi_a(t)|$. The smoothing equation is
\[
 \phi_a(t)=\int_0^1\phi_a(tv^a)
 \{p\phi_a(t(1-v)^a)+(1-p)\phi_a(-t(1-v)^a)\}\,dv,
 \qquad p=(1+a)/2.
\]
The local bounds just proved justify differentiation under this integral,
since $v^{\min J}|\log v|$ is integrable. Put $x=tv^a$,
$y=t(1-v)^a$ and, for nonnegative $z$, define
\[
 (\mathcal K_a z)(t)=\int_0^1\{z(x)f_a(y)+f_a(x)z(y)\}\,dv.
\]
Taking absolute values and discarding factors $v^a,(1-v)^a\le1$
in the frequency derivative gives
\[
 u_a\le\mathcal K_a u_a,\qquad
 w_a\le\mathcal K_a w_a+S_a,
\]
where
\[
 S_a(t)=\int_0^1\{x|\log v|u_a(x)f_a(y)
       +y|\log(1-v)|f_a(x)u_a(y)+f_a(x)f_a(y)\}\,dv.
\]
The last term bounds the derivative of $p$, since $p'=1/2$ and
$|\phi_a(y)-\phi_a(-y)|\le2f_a(y)$.

Fix $\beta>0$. Splitting at $1/2$ and reflecting the second half gives
\[
 K_{a,\beta}(t)=2\int_0^{1/2}
 \left\{\left(\frac{1+t}{1+x}\right)^\beta f_a(y)
       +\left(\frac{1+t}{1+y}\right)^\beta f_a(x)\right\}\,dv
 \longrightarrow0
\]
uniformly for $a\in J$ as $t\to\infty$.
Indeed, $y\ge t/2$ bounds the first term by
$2(1+t)^\beta e^{-c(t/2)^r}$ and the second ratio by $2^\beta$.
Also $x\ge tv^{\max J}$, so the exponential bound and dominated
convergence make the second integral tend to zero uniformly.
Choose $T$ so that $K_{a,\beta}(t)\le1/2$ for all $a\in J$, $t\ge T$.
For $X\ge T$ let
\[
 U_X=\max_{a\in J,\,0\le t\le X}(1+t)^\beta u_a(t).
\]
All arguments $x,y$ are at most $t$. The first smoothing inequality gives
\[
 U_X\le\max\{U_T,U_X/2\},\qquad\text{hence }U_X\le U_T.
\]
The bound is independent of $X$. This proves arbitrary polynomial decay for $u_a$.

The same decay holds for $S_a$. On $0<v\le1/2$, the first
logarithmic term uses the exponential decay of $f_a(y)$ and boundedness
of $u_a(x)$, with $x|\log v|\le tv^{\min J}|\log v|$.
The second uses arbitrary polynomial decay of $u_a(y)$,
$f_a(x)\le1$, and $|\log(1-v)|\le\log2$. The last uses the
exponential decay of $f_a(y)$. Reflection covers the other half.
Thus $S_\beta:=\sup_{a\in J,t\ge0}(1+t)^\beta S_a(t)<\infty$.
Defining $W_X$ as above with $w_a$ in place of $u_a$, we obtain
\[
 W_X\le\max\{W_T,W_X/2+S_\beta\}
 \quad\Longrightarrow\quad W_X\le\max\{W_T,2S_\beta\}.
\]
Conjugate symmetry extends both derivative bounds to negative $t$.
A finite cover gives the stated bounds on every compact memory interval.
Since $\beta$ is arbitrary, multiplication by any fixed power of $t$
leaves an integrable uniform majorant. Differentiation under Fourier
inversion and dominated convergence prove the claimed joint continuity
of $\partial_a\partial_u^k h_a$ for every $k$.
\end{proof}

\section{The quadratic rate near one}\label{app:quadratic-rate}

\begin{proof}[Proof of Proposition~\ref{prop:quadratic-rate}]
We compare the mass of $g_{a,q}$ on two adjacent intervals $J_1$
and $J_2$ of equal length to the left of zero. The reflected law puts
mass at least $\varpi$ on $J_1$, and the direct law puts mass at most
$A$ on $J_2$. The number $A$ is the product of two factors of order
$1-a$, the relative width of a window and the probability that
$\mathcal L_a\le0$. This is why the rate is quadratic.
Write $\Prob_+$ for the law of the walk given $X_1=1$. Let
$p=1-\varepsilon$, fix an integer $k\ge1$, and let
$E_k=\{X_1=\dots=X_k=1\}$. Given $X_1=1$, each of the next $k-1$ steps
copies a positive step with probability $1-\varepsilon$ or flips a
positive step with probability $\varepsilon$, so
$\Prob_+(E_k)=p^{k-1}$. After time $k$, the step $n+1$ is $+1$ with
probability $\frac12+\frac{aS_n}{2n}$. Hence
\[
 \E[S_{n+1}\mid\mathcal F_n]=\Bigl(1+\frac an\Bigr)S_n,\qquad
 \E[S_{n+1}^2\mid\mathcal F_n]=\Bigl(1+\frac{2a}n\Bigr)S_n^2+1.
\]
Put $\gamma_n=\prod_{j=k}^{n-1}(1+a/j)$ and
$\beta_n=\prod_{j=k}^{n-1}(1+2a/j)$. Given $E_k$, the first relation
makes $S_n/\gamma_n$ a martingale with mean $k$, and the second gives
$\E[S_n^2\mid E_k]=\beta_n\bigl(k^2+\sum_{m=k+1}^n\beta_m^{-1}\bigr)$.
The identity
\[
 (2a-1)\frac{\Gamma(m)}{\Gamma(m+2a)}
 =\frac{\Gamma(m)}{\Gamma(m+2a-1)}-\frac{\Gamma(m+1)}{\Gamma(m+2a)}
\]
telescopes to $\sum_{m>k}\beta_m^{-1}=k/(2a-1)$. Since
$\gamma_n/n^a\to\Gamma(k)/\Gamma(k+a)$ and
$\beta_n/\gamma_n^2\to c_k=\Gamma(k+a)^2/(\Gamma(k)\Gamma(k+2a))$, the
martingale is bounded in $L^2$. Its limit gives, given $E_k$,
\[
 \E[\mathcal L_a\mid E_k]=m_k=\frac{\Gamma(k+1)}{\Gamma(k+a)},\qquad
 \Var(\mathcal L_a\mid E_k)=r_km_k^2,\qquad
 r_k=c_k\Bigl(1+\frac1{(2a-1)k}\Bigr)-1.
\]
For $k=1$ this is $r_1=\rho_a$. To bound $r_k$, let $\varphi(s)$ be
$\log(1+r_k)$ with $a$ replaced by $s$. Then $\varphi(1)=0$ and
\[
 \varphi'(s)=-2\bigl(\psi(k+2s)-\psi(k+s)\bigr)
 -\frac2{(2s-1)\{(2s-1)k+1\}}<0.
\]
Since $\psi'(x)<1/(x-1)$ for $x>1$, we have
$\psi(k+2s)-\psi(k+s)\le s\psi'(k+s)\le1/k$ for $s\le1$. Also
$(2s-1)k+1\ge(2s-1)(k+1)$. So $0\le\log(1+r_k)\le(1-a)\theta_k$ with
$\theta_k=2/k+2/((2a-1)^2(k+1))$.

Now fix $K>1$, put $v=\sqrt{r_k}$, and assume $3Kv<1$ and
$\varpi=p^{k-1}(1-K^{-2})>1/2$. Let $x_0=m_k(1-Kv)$ and $\ell=2Km_kv$, so
$x_0-\ell=m_k(1-3Kv)>0$. By Chebyshev's inequality given $E_k$, the
window $W=[x_0,x_0+\ell]$ has $\Prob_+(\mathcal L_a\in W)\ge\varpi$, and
so $\Prob_+(\mathcal L_a<x_0)\le1-\varpi$. By
Lemmas~\ref{lem:balanced-start} and~\ref{lem:positive-conditional-mode},
$h_a$ is increasing on $(-\infty,0]$. So for $y<0$ we have
$|y|h_a(y)\le\Prob_+(\mathcal L_a\le0)$, and Cantelli's inequality
gives $\Prob_+(\mathcal L_a\le0)\le\rho_a/(1+\rho_a)$. On
$J_2=[-x_0,-x_0+\ell]$ this gives
\[
 \Prob_+(\mathcal L_a\in J_2)\le\frac{\ell}{x_0-\ell}\cdot
 \frac{\rho_a}{1+\rho_a}
 =\frac{2Kv}{1-3Kv}\cdot\frac{\rho_a}{1+\rho_a}=:A.
\]
Let $q\in[1/2,1)$ and $J_1=[-x_0-\ell,-x_0]$, so $-J_1=W$ and
$-J_2\subset(-\infty,x_0]$. Then
\begin{align*}
 \int_{J_1}g_{a,q}-\int_{J_2}g_{a,q}
 &\ge-q\,\Prob_+(\mathcal L_a\in J_2)
   +(1-q)\bigl(\Prob_+(\mathcal L_a\in W)-\Prob_+(\mathcal L_a<x_0)\bigr)\\
 &\ge-qA+(1-q)(2\varpi-1).
\end{align*}
If $A\ge2\varpi-1$, the bound below holds because $b_a\le1/2$. So let
$A<2\varpi-1$. If $1-q>A/(2\varpi-1+A)$, the right side is positive. The intervals $J_1$
and $J_2$ are adjacent, have equal length, and lie in $(-\infty,0)$
with $J_1$ on the left, so $g_{a,q}$ is not nondecreasing on
$(-\infty,0]$. Thus $g_{a,q}'(-t)=qh_a'(-t)-(1-q)h_a'(t)<0$ for some
$t>0$. Since $h_a'(-t)>0$, this forces $t\in U_a$ and
$h_a'(-t)/(h_a'(t)+h_a'(-t))<1-q$. By~\eqref{eq:bias-boundary},
$b_a<1-q$, and therefore
\[
 b_a\le\frac{A}{2\varpi-1+A}\le\frac A{2\varpi-1}.
\]

Finally let $a\uparrow1$. A run of $k$ positive first steps with
$k$ of order $1/(1-a)$ has probability bounded below, and it lowers
the relative variance from order $1-a$ to order $(1-a)^2$. The
constants below were chosen to make the final bound small, and we do
not claim that $82$ is optimal. Take $\delta=4/25$, $K=27/10$, and
$k=\lceil2\delta/(1-a)\rceil$. Then $p^{k-1}\to e^{-\delta}$ and
$k\theta_k\to4$, so $v^2\le e^{(1-a)\theta_k}-1=(2/\delta+o(1))(1-a)^2$.
With $\rho_a=2(1-a)+O((1-a)^2)$ this gives
$A=(4K\sqrt{2/\delta}+o(1))(1-a)^2$ and
$2\varpi-1\to2e^{-\delta}(1-K^{-2})-1>0$. Hence
\[
 \limsup_{a\uparrow1}\frac{b_a}{(1-a)^2}
 \le\frac{4K\sqrt{2/\delta}}{2e^{-\delta}(1-K^{-2})-1}<82.
\]
\end{proof}

\section*{Statement on AI use}
I used two AI tools for this paper, with their contributions described here. Claude Opus 5.5
(Anthropic, through Claude Code) ran literature searches, citation checks,
and referee-style reviews, checked computations in the threshold and
quartic-top proofs, and made final copy edits and link updates. The
reviews informed the framing and the shortening of an earlier combined
manuscript, and suggested proving without computation that the bias
boundary is positive, that the transition is unique, and that the critical
top is quartic. Later reviews of the exposition led me to reorder the
paper, add plans to the long proofs, rename several symbols, and they pointed out the gaps that I closed in the proofs of
Lemmas~A.3 and~A.4 and of the proposition in Appendix~E. All proofs in
the paper were written by me, and the aforementioned suggestions are the
only Claude suggestions that led to a new proof or result. I wrote the
Python scripts in \nolinkurl{certificates}, and Claude Opus 5.5 revised
them and redrew Figure~1 to be more descriptive. The Lean files and their coverage map in
\nolinkurl{lean/PaperD} were generated with Claude Opus 5.5 from my
statements and proofs, and I checked that they build and audited their
axioms. GPT-6 Astra (OpenAI, through Codex) edited my
drafts for grammar, concision, and consistent notation, checked the
citations, and read the proofs for errors. It reported inconsistent
notation, and it helped me close a gap, after the setup, in my first proof
that $b_a>0$. Neither tool wrote a section of the paper. I read every cited
source in full except Eisenberger (1964) and Robertson and Fryer (1969),
which are cited only for results stated in their abstracts.
Section~\ref{sec:verify} separates the written proofs, the formal proofs
with their named classical hypotheses, and the numerical certificates.
I take full responsibility for the paper, including its proofs and
computations.


\begin{thebibliography}{99}\small\setlength{\baselineskip}{11pt}\setlength{\itemsep}{0pt}\setlength{\parsep}{0pt}

\bibitem{schutz2004}
G.~M.~Sch\"utz and S.~Trimper, \emph{Elephants can always remember: exact
long-range memory effects in a non-Markovian random walk}, Phys.\ Rev.\ E
\textbf{70} (2004), 045101(R).
\url{https://doi.org/10.1103/PhysRevE.70.045101}.

\bibitem{dreznerfarnum1993}
Z.~Drezner and N.~Farnum, \emph{A generalized binomial distribution},
Comm.\ Statist.\ Theory Methods \textbf{22} (1993), 3051--3063.
\url{https://doi.org/10.1080/03610929308831202}.

\bibitem{heyde2004}
C.~C.~Heyde, \emph{Asymptotics and criticality for a correlated Bernoulli
process}, Aust.\ N.\ Z.\ J.\ Statist.\ \textbf{46} (2004), 53--57.
\url{https://doi.org/10.1111/j.1467-842X.2004.00311.x}.

\bibitem{jamesjamesqi2008}
B.~James, K.~James and Y.~Qi, \emph{Limit theorems for correlated Bernoulli
random variables}, Statist.\ Probab.\ Lett.\ \textbf{78} (2008), 2339--2345.
\url{https://doi.org/10.1016/j.spl.2008.01.104}.

\bibitem{coletti2017}
C.~F.~Coletti, R.~Gava and G.~M.~Sch\"utz, \emph{Central limit theorem and
related results for the elephant random walk}, J.~Math.~Phys. \textbf{58}
(2017), 053303.
\url{https://doi.org/10.1063/1.4983566}.

\bibitem{baurbert2016}
E.~Baur and J.~Bertoin, \emph{Elephant random walks and their
connection to P\'olya-type urns}, Phys.\ Rev.\ E \textbf{94}
(2016), 052134. \url{https://doi.org/10.1103/PhysRevE.94.052134}.

\bibitem{kursten2016}
R.~K\"ursten, \emph{Random recursive trees and the elephant random walk},
Phys.\ Rev.\ E \textbf{93} (2016), 032111.
\url{https://doi.org/10.1103/PhysRevE.93.032111}.

\bibitem{bertoin2022}
J.~Bertoin, \emph{Counting the zeros of an elephant random walk},
Trans.\ Amer.\ Math.\ Soc.\ \textbf{375} (2022), 5539--5560.
\url{https://doi.org/10.1090/tran/8622}.

\bibitem{qin2025}
S.~Qin, \emph{Recurrence and transience of multidimensional elephant
random walks}, Ann.\ Probab.\ \textbf{53} (2025), 1049--1078.
\url{https://doi.org/10.1214/24-AOP1728}.

\bibitem{bercu2018}
B.~Bercu, \emph{A martingale approach for the elephant random walk},
J.\ Phys.\ A: Math.\ Theor.\ \textbf{51} (2018), 015201.
\url{https://doi.org/10.1088/1751-8121/aa95a6}.

\bibitem{kubotatakei2019}
N.~Kubota and M.~Takei, \emph{Gaussian fluctuation for superdiffusive
elephant random walks}, J.~Stat.~Phys. \textbf{177} (2019), 1157--1171.
\url{https://doi.org/10.1007/s10955-019-02414-0}.

\bibitem{glrs2025}
H.~Gu\'erin, L.~Laulin, K.~Raschel and T.~Simon, \emph{On the limit law of
the superdiffusive elephant random walk}, Electron.\ J.\ Probab.\
\textbf{30} (2025), paper no.~102, 1--25. \url{https://doi.org/10.1214/25-EJP1362}.

\bibitem{glr2026}
H.~Gu\'erin, L.~Laulin and K.~Raschel, \emph{A fixed-point equation approach
for the superdiffusive elephant random walk}, Ann.\ Inst.\ Henri Poincar\'e
Probab.\ Stat.\ \textbf{62} (2026), 973--1005.
\url{https://doi.org/10.1214/24-AIHP1534}.

\bibitem{kissveto2022}
J.~Kiss and B.~Vet\H{o}, \emph{Moments of the superdiffusive elephant
random walk with general step distribution}, Electron.\ Commun.\ Probab.\
\textbf{27} (2022), paper no.~44, 1--12.
\url{https://doi.org/10.1214/22-ECP485}.

\bibitem{peresqin2026}
Y.~Peres and S.~Qin, \emph{Transition probabilities of step-reinforced
random walks}, preprint, 2026.
\href{https://arxiv.org/abs/2604.07227}{arXiv:2604.07227}.

\bibitem{tv2026}
H.~Trinas and G.~Valle, \emph{A superdiffusive local limit theorem for the
elephant random walk and the breakdown of log-concavity}, preprint, 2026.
\href{https://arxiv.org/abs/2609.14753}{arXiv:2609.14753}.

\bibitem{trinas2026}
H.~Trinas, \emph{The prodigal elephant returns: last-return laws in the
superdiffusive regime}, preprint, 2026.
\href{https://arxiv.org/abs/2611.01279}{arXiv:2611.01279}.

\bibitem{muraharatokumitsu2026}
H.~Murahara and G.~Tokumitsu, \emph{Optimal rates in local and central limit
theorems for the elephant random walk}, preprint, 2026.
\href{https://arxiv.org/abs/2610.35019}{arXiv:2610.35019}.

\bibitem{vishwakarma2026}
P.~Vishwakarma, \emph{Moments approach for the elephant random walk},
preprint, 2026.
\href{https://arxiv.org/abs/2609.11035}{arXiv:2609.11035}.

\bibitem{kagey131}
P.~Kagey, \emph{Problem 131}, Open Problem Collection,
\url{https://peterkagey.com/problems/131/} (accessed October 2026).

\bibitem{johansson2017}
F.~Johansson, \emph{Arb: efficient arbitrary-precision midpoint-radius
interval arithmetic}, IEEE Trans.\ Comput.\ \textbf{66} (2017), 1281--1292.
\url{https://doi.org/10.1109/TC.2017.2690633}.

\bibitem{dj1988}
S.~Dharmadhikari and K.~Joag-dev, \emph{Unimodality, Convexity, and
Applications}, Academic Press, Boston, 1988.

\bibitem{dasilva2013}
M.~A.~A.~da~Silva, J.~C.~Cressoni, G.~M.~Sch\"utz, G.~M.~Viswanathan and
S.~Trimper, \emph{Non-Gaussian propagator for elephant random walks},
Phys.\ Rev.\ E \textbf{88} (2013), 022115.
\url{https://doi.org/10.1103/PhysRevE.88.022115}.

\bibitem{mohandas2025}
K.~Mohandas, P.~J.~G\'orski, K.~Suchecki, G.~Andres, G.~Vaccario and J.~A.~Ho{\l}yst,
\emph{Structurally balanced growing network as randomized P\'olya urn process},
preprint, 2025.
\href{https://arxiv.org/abs/2510.24659}{arXiv:2510.24659}.

\bibitem{karlin1968}
S.~Karlin, \emph{Total Positivity}, Vol.~I, Stanford University Press,
Stanford, CA, 1968.

\bibitem{karlin1964}
S.~Karlin, \emph{Total positivity, absorption probabilities and applications},
Trans.\ Amer.\ Math.\ Soc.\ \textbf{111} (1964), 33--107.
\url{https://doi.org/10.1090/S0002-9947-1964-0168010-2}.

\bibitem{eisenberger1964}
I.~Eisenberger, \emph{Genesis of bimodal distributions},
Technometrics \textbf{6} (1964), 357--363.
\url{https://doi.org/10.1080/00401706.1964.10490199}.

\bibitem{robertsonfryer1969}
C.~A.~Robertson and J.~G.~Fryer, \emph{Some descriptive properties of
normal mixtures}, Scand.\ Actuar.\ J.\ \textbf{1969} (1969), 137--146.
\url{https://doi.org/10.1080/03461238.1969.10404590}.
\end{thebibliography}
\end{document}